\documentclass[11pt]{article}

\usepackage{amsmath}
\usepackage{amssymb}
\usepackage{amsfonts}
\usepackage{amsthm}
\usepackage{amscd}

\usepackage{graphicx}
\usepackage{color}
\usepackage{xcolor}
\usepackage{tikz}
\usepackage{pgfplots}
\pgfplotsset{width=10cm,compat=1.9}
\usepackage{hyperref}

\usepackage{hyperref}
\usepackage{caption}
\usepackage{subcaption}
\usepackage{array}
\usepackage{bbm}
\usepackage{listings}
\usepackage{complexity}
\usepackage[ruled,vlined]{algorithm2e}
\usepackage{lineno}
\usepackage[margin=0.25in]{geometry}

\newtheorem{thm}{Theorem}
\newtheorem{conj}[thm]{Conjecture}

\newtheorem{lem}[thm]{Lemma}
\newtheorem{quest}{Question}

\newtheorem{swap}{Extension Rule}

\newtheorem{claim}{Claim}
\newtheorem{subclaim}{Claim}[claim]

\newcommand{\smallqed}{{\tiny ($\Box$)}}

\newcommand{\QEDmark}{\mbox{\textsc{qed}}}
\newcommand{\proofStarter}[1]{\textsc{#1} }

\def\vertex(#1){\put(#1){\circle*{2}}}
\def\vertexo(#1){\put(#1){\circle{2}}}
\def\vert(#1){\put(#1){\circle*{1.5}}}
\def\verto(#1){\put(#1){\circle{1.5}}}
\def\lab(#1)#2{\put(#1){\makebox(0,0)[c]{#2}}}
\definecolor{DarkGreen}{rgb}{0.2, 0.6, 0.3}

\definecolor{electricindigo}{rgb}{0.44, 0.0, 1.0}

\let\oldenumerate\enumerate
\renewcommand{\enumerate}{
  \oldenumerate
  \setlength{\itemsep}{0.5pt}
  \setlength{\parskip}{0pt}
  \setlength{\parsep}{0pt}
}

\begin{document}

\title{Another conjecture of \emph{TxGraffiti} concerning \\ zero forcing and domination in graphs}
\author{$^{1,2}$Randy Davila\
\\
$^1$Research and Development \\
RelationalAI \\
Berkeley, CA 94704, USA\\
\small {\tt Email: randy.davila@relational.ai}\\
\\
$^2$Department of Computational Applied \\ Mathematics \& Operations Research\\
Rice University\\
Houston, TX 77005, USA \\
\small {\tt Email: randy.r.davila@rice.edu} \\
\\
}

\date{}
\maketitle

\begin{abstract}
This paper proves a conjecture generated by the artificial intelligence conjecturing program called \emph{TxGraffiti}. More specifically, we show that if $G$ is a connected, cubic, and claw-free graph, then $Z(G) \le \gamma(G) + 2$, where $Z(G)$ and $\gamma(G)$ denote the zero forcing number and the domination number of $G$, respectively. Furthermore, we provide a complete characterization of graphs that achieve this bound. Notably, this bound improves the known upper bounds for the zero forcing number of connected, cubic, and claw-free graphs.
\end{abstract}


\section{Introduction}
Written by the author in 2017 and tailored towards graph theoretic conjectures, \emph{TxGraffiti} is a machine learning and artificially intelligent program for generating mathematical conjectures suitable for research by professional mathematicians; now available as an interactive website\footnote{\href{https://txgraffiti.streamlit.app}{https://txgraffiti.streamlit.app}}. The name \emph{TxGraffiti} plays a homage to both the original conjecturing program \emph{Graffiti} written by Fajtlowicz~\cite{Graffiti}, and its successor program \emph{Graffiti.pc} written by DeLaViña~\cite{Graffiti.pc}. However, the design of \emph{TxGraffiti} remains distinct from both \emph{Graffiti} and \emph{Graffiti.pc}, so the conjectures of \emph{TxGraffiti} differ in many instances. One particularly interesting difference is the ability of \emph{TxGraffiti} to make particularly strong and interesting conjectures on \emph{domination} and \emph{zero forcing}-related parameters for cubic graphs. One such example is that if $G \neq K_4$ is a connected and cubic graph, then $Z(G) \leq 2\gamma(G)$, where $Z(G)$ and $\gamma(G)$ denote the \emph{zero forcing number} and \emph{domination number} of $G$, respectively, a conjecture that was confirmed in~\cite{Davila2, DaHe21a}. Notably, a similar conjecture of \emph{TxGraffiti} relating \emph{total zero forcing} and \emph{total domination} for cubic graphs was also shown to be true in~\cite{Davila2, DaHe19b}. A more recent conjecture on zero-forcing made by \emph{TxGraffiti} states that the zero forcing number is, at most, the vertex cover number for claw-free graphs, which in turn was confirmed and led to a larger body of work presented by Brimkov et al.~\cite{BrDaScYo}. For \emph{TxGraffiti} inspired results that do not concern the zero forcing number of a graph, see~\cite{CaDaHePe2022}, or see the surprising result given by Caro et al.~\cite{CaDaPe}, which states that $\alpha(G)\leq \mu(G)$ for any regular graph $G$, where $\alpha(G)$ and $\mu(G)$ denote the \emph{independence number} and \emph{matching number} of $G$, respectively. 

Though many conjectures of \emph{TxGraffiti} are confirmed, many more remain open, and even more remain \emph{undiscovered}. The most well-known of these open conjectures is likely \emph{the $\alpha$-Z Conjecture}; namely, if $G\neq K_4$ is a connected and cubic graph, then $Z(G) \leq \alpha(G) + 1$; see the last page of Davila's dissertation~\cite{Davila2}. The $\alpha$-Z conjecture has received a considerable amount of attention with several partial results, including $Z(G) \leq \alpha(G) + 1$, whenever $G \neq K_4$ is connected, cubic, and claw-free~\cite{Davila2, DaHe19c}. Notably, this partial result was \emph{not conjectured} by \emph{TxGraffiti}, and so, by the design of \emph{TxGraffiti}, there must exist a more substantial upper bound for $Z(G)$ in cubic and claw-free graphs in terms of a graph parameter other than $\alpha(G)$. Thus, we discovered the following stronger conjecture when explicitly examining the conjectures of \emph{TxGraffiti} for connected, cubic, and claw-free graphs. 
\begin{conj}[\emph{TxGraffiti} -- confirmed]\label{conj:confirmed}
If $G$ is a connected, cubic, and claw-free graph, then
\[
Z(G) \leq \gamma(G) + 2, 
\]
and this bound is sharp. 
\end{conj}

In this paper, we resolve Conjecture~\ref{conj:confirmed} in the affirmative, and in doing so, also characterize graphs attaining equality in its statement. Notably, the resulting theorem improves on both $Z(G) \leq \alpha(G) + 1$ and $Z(G) \leq 2\gamma(G)$ (for $G \neq K_4$) whenever $G$ is a connected, cubic, and claw-free graph. 

\subsection{Notation and Terminology}
Throughout this paper, all graphs considered are simple, undirected, and finite. Let $G$ be a graph with vertex set $V(G)$ and edge set $E(G)$. The \emph{order} of $G$ is $n(G) = |V(G)|$. Two vertices $v,w \in V(G)$ are \emph{neighbors}, or \emph{adjacent}, if $vw \in E(G)$. The \emph{open neighborhood} of $v\in V(G)$, is the set of neighbors of $v$, denoted $N_G(v)$. The closed neighborhood of $S\subseteq V$ is $N_G[S] = \bigcup_{v\in S}N_G[v]$. The \emph{degree} of a vertex $v\in V(G)$, denoted $d_G(v)$, is equal to $|N_G(v)|$. A \emph{cubic graph} (also called a \emph{3-regular graph}) is a graph for which every vertex degree is three. The maximum and minimum degree of $G$ will be denoted $\Delta(G)$ and $\delta(G)$, respectively. When there is no scope for confusion, we will use the notation $n = n(G)$, $\delta = \delta(G)$, and $\Delta = \Delta(G)$, to denote the order, minimum degree, and maximum degree, respectively. 

Two vertices in a graph $G$ are \emph{independent} if they are not neighbors. A set of pairwise independent vertices in $G$ is an \emph{independent set} of $G$. The number of vertices in a maximum independent set in $G$ is the \emph{independence number} of $G$, denoted $\alpha(G)$. A \emph{dominating set} in $G$ is a subset $X\subseteq V(G)$, so all vertices in $G$ are either in $X$ or adjacent with at least one vertex in $X$. The \emph{domination number} of $G$, denoted by $\gamma(G)$, is the cardinality of a minimum dominating set in $G$. If $X \subseteq V(G)$ is a dominating and independent set, then $X$ is an \emph{independent dominating set}. The minimum cardinality of an independent dominating set in $G$ is the \emph{independent domination number} of $G$, denoted $i(G)$. For more on domination in graphs, see the comprehensive and well-written text by Haynes, Hedetniemi, and Henning~\cite{HaHeHe2024}. 

For a set of vertices $S \subseteq V(G)$, the subgraph induced by $S$ is denoted by $G[S]$. If $v\in V(G)$, we denote the graph obtained by deleting $v$ in $G$ by $G-v$. We denote the path, cycle, and complete graph on $n$ vertices by $P_n$, $C_n$, and $K_4$, respectively. A \emph{triangle} in $G$ is an induced subgraph of $G$ isomorphic to $K_3$, whereas a \emph{diamond} in $G$ is a subgraph of $G$ isomorphic to $K_4$ with one edge missing. A graph $G$ is \emph{$F$-free} if $G$ does not contain $F$ as an induced subgraph. In particular, if $G$ is $F$-free, where $F = K_{1,3}$, then $G$ is \emph{claw-free}. Claw-free graphs have been widely studied, and the survey by Flandrin, Faudree, and Ryj{\'a}{\v{c}}ek~\cite{claw-free} is the standard reference for the topic. 

Let $S\subseteq V(G)$ be a set of initially ``blue-colored'' vertices, all remaining vertices being ``white-colored''. The \emph{zero forcing process} on $G$ is defined as follows: At each discrete time step, if a blue-colored vertex has a unique white-colored neighbor, then this blue-colored vertex \emph{forces} its white-colored neighbor to become colored blue. The \emph{open $z$-neighborhood} of $S$ in $G$, written $N_G^z(S)$, is the set of all initially white-colored vertices that change color during the zero forcing process. The \emph{closed $z$-neighborhood} of $S$ in $G$ is the set $N_G^z[S] = S \cup N_G^z(S)$. A \emph{valid-$z$-extension rule} on $S$ is a mapping from $S$ to a new set $S'$, such that $|S| = |S'|$ and $N_G^z[S] \subseteq N_G^z[S']$; if in addition $N_G^z[S] \subset N_G^z[S']$, then we say the rule is a \emph{valid-$z$-proper-extension rule} on $S$. A set $X \subseteq V(G)$ is \emph{$S$-reachable} when $X\subseteq N_G^z[S]$. If $V(G)$ is $S$-reachable, and so, $V(G) = N_G^z[S]$, then $S$ is a \emph{zero forcing set} of $G$. The minimum cardinality of a zero forcing set in $G$ is the \emph{zero forcing number} of $G$. Zero forcing was introduced in \cite{AIM-Workshop} in the context of linear algebra and has since gained much interest in graph theory. For more information on zero forcing in graphs, we refer the reader to the excellent monograph by Hogben, Lin, and Shader~\cite{HoLiSh-zero-forcing-book}. 

For notation and graph terminology not mentioned here, we refer the reader to~\cite{HaHeHe2024}. We also will use the standard notation $[k] =\{1,\ldots,k\}$. 


\section{Known Results and Lemmas}\label{sec:lemmas}
In this section, we introduce terminology, recall useful results, and prove technical lemmas needed for the proof of Conjecture~\ref{conj:confirmed}. To begin, recall the following result of Allan and Laskar~\cite{AlLa1978}, which states that the independent domination number equals the domination number for all claw-free graphs. 
\begin{thm}[\cite{AlLa1978}]\label{thm:i-equal-gamma}
If $G$ is a claw-free graph, then $i(G) = \gamma(G)$. 
\end{thm}

We will also need the following structural property of connected, cubic, and claw-free graphs established by Henning and L\"{o}wenstein in~\cite{HeLo12}.
\begin{lem}[\cite{HeLo12}]\label{lem:partition}
If $G\neq K_4$ is a connected, cubic, and claw-free graph, then the vertex set $V(G)$ can be uniquely partitioned into sets, each of which induces a triangle or diamond in $G$.
\end{lem}

By Lemma~\ref{lem:partition}, the vertex set $V(G)$ of a connected, cubic, and claw-free graph $G\neq K_4$ can be uniquely partitioned into sets, each of which induces either a triangle subgraph or a diamond subgraph of $G$. We refer to such a partition of $G$ as a \emph{triangle-diamond} partition of $G$, abbreviated $\Delta$-D-partition. Borrowing the terminology introduced in~\cite{HeLo12}, we call every triangle and diamond induced by our $\Delta$-D-partition a \emph{unit} of the partition. A unit that is a triangle is called a \emph{triangle-unit} and a unit that is a diamond is called a \emph{diamond-unit}. Note that triangle-units do not belong to any diamond-unit. Furthermore, we say that two units are adjacent if an edge joins a vertex in one unit to a vertex in another and that a vertex is adjacent to a unit whenever the vertex has a neighbor. 

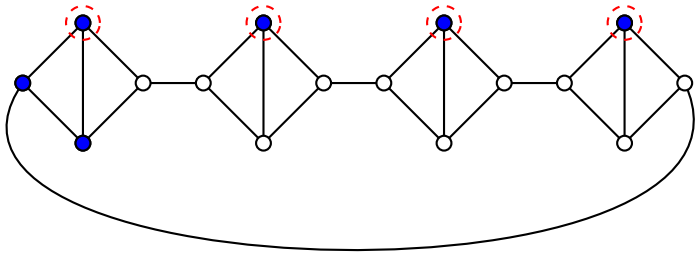
\begin{figure}[htb]
\begin{center}
\begin{tikzpicture}[scale=.8,style=thick,x=1cm,y=1cm]
\def\vr{3.5pt}

\path (-.5, 0.0) coordinate (d11);
\path (0.5, 1.0) coordinate (d12);
\path (1.5, 0.0) coordinate (d13);
\path (0.5, -1.0) coordinate (d14);

\path (2.5, 0.0) coordinate (d21);
\path (3.5, 1.0) coordinate (d22);
\path (4.5, 0.0) coordinate (d23);
\path (3.5, -1.0) coordinate (d24);

\path (5.5, 0.0) coordinate (d31);
\path (6.5, 1.0) coordinate (d32);
\path (7.5, 0.0) coordinate (d33);
\path (6.5, -1.0) coordinate (d34);

\path (8.5, 0.0) coordinate (d41);
\path (9.5, 1.0) coordinate (d42);
\path (10.5,0.0) coordinate (d43);
\path (9.5, -1.0) coordinate (d44);

\draw (d11) -- (d12);
\draw (d12) -- (d13);
\draw (d11) -- (d14);
\draw (d14) -- (d13);
\draw (d12) -- (d14);

\draw (d13) -- (d21);

\draw (d21) -- (d22);
\draw (d22) -- (d23);
\draw (d21) -- (d24);
\draw (d24) -- (d23);
\draw (d22) -- (d24);

\draw (d23) -- (d31);

\draw (d31) -- (d32);
\draw (d32) -- (d33);
\draw (d31) -- (d34);
\draw (d34) -- (d33);
\draw (d32) -- (d34);

\draw (d33) -- (d41);

\draw (d41) -- (d42);
\draw (d42) -- (d43);
\draw (d41) -- (d44);
\draw (d44) -- (d43);
\draw (d42) -- (d44);

\draw (d11) to[out=-125,in=-65] (d43);

\foreach \i in {1,2,3,4} {
    \draw (d\i1) [fill=white] circle (\vr);
    \draw (d\i2) [fill=white] circle (\vr);
    \draw (d\i3) [fill=white] circle (\vr);
    \draw (d\i4) [fill=white] circle (\vr);
}

\draw (d12) [fill=blue] circle (\vr);
\draw[dashed, red] (d12) circle (8pt);

\draw (d22) [fill=blue] circle (\vr);
\draw[dashed, red] (d22) circle (8pt);

\draw (d32) [fill=blue] circle (\vr);
\draw[dashed, red] (d32) circle (8pt);

\draw (d42) [fill=blue] circle (\vr);
\draw[dashed, red] (d42) circle (8pt);

\draw (d11) [fill=blue] circle (\vr);
\draw (d14) [fill=blue] circle (\vr);

\end{tikzpicture}
\end{center}
\caption{The diamond-necklace $N_4$ with a minimum zero forcing set shown by blue colored vertices and a minimum dominating set shown by the red dashed circles.}
\label{fig:diamond-necklace}
\end{figure}

Henning and L\"{o}wenstein also introduced the notation and definition of a \emph{diamond-necklace}, which is now restated. For $k \geq 2$ an integer, let $N_k$ be the connected and cubic graph constructed as follows. Take $k$ disjoint copies $D_1, D_2, . . . , D_k$ of 
a diamond, where $V(D_i) = \{a_i, b_i, c_i, d_i \}$ and where $a_ib_i$ is the missing edge in $D_i$. Let $N_k$ be obtained from the disjoint union of 
these $k$ diamonds by adding the edges $\{a_{i}b_{i+1} \: | \: i \in [k-1]\}$ and adding the edge $a_kb_1$. We call $N_k$ a \emph{diamond-necklace} with $k$ diamonds; see Figure~\ref{fig:diamond-necklace} for the diamond-necklace $N_4$. Let $\mathcal{N}_{cubic} = \{N_k \: | \:  k \geq 2\}$. The complete graph $K_4$ may actually be thought of as a loop diamond-necklace, and so, we define $\mathcal{N}_{cubic}^* = \mathcal{N}_{cubic} \cup \{K_4\}$. In~\cite{Davila2, DaHe18b}, it was established that if $G \in \mathcal{N}_{cubic}$, then $Z(G) = \frac{1}{4}n(G) + 2$. This together with the fact that $Z(K_4) = 3 = \frac{1}{4}n(K_4) + 2$, imply the following lemma. 
\begin{lem}\label{lem:order-bound}
If $G \in \mathcal{N}_{cubic}^*$ has order $n$, then $Z(G) = \frac{1}{4}n + 2$. 
\end{lem}

Let $D$ be a diamond-unit in the $\Delta$-D-partition of a cubic and claw-free graph $G$ and let $V(D) = \{a, b, c, d\}$ denote the vertex set of $D$, where $ab$ is the missing edge in $D$. Then, the only neighbors of vertices $b$ and $c$, are vertices contained in $V(D)$, and as a consequence of this, we obtain the following lemma. 
\begin{lem}\label{lem:diamond-containment}
If $D$ is a diamond-unit in a $\Delta$-D-partition of the cubic and claw-free graph $G$, and if $X \subseteq V(G)$ is a dominating set of $G$, then $V(D) \cap X \neq \emptyset$. 
\end{lem}

A trivial lower bound for the domination number of any graph $G$ of order $n$ and maximum degree $\Delta$, is $\gamma(G) \geq \frac{n}{\Delta+1}$. Thus, for a cubic graph $G$ of order $n$, we have the lower bound $\gamma(G) \geq \frac{n}{4}$. Let $D_1, D_2, . . . , D_k$ be the diamond-units in the diamond-necklace $N_k$, where $V(D_i) = \{a_i, b_i, c_i, d_i \}$ and where $a_ib_i$ is the missing edge in $D_i$, and let $n$ be the order of $N_k$. It is easy to see that $X = \{c_1, \dots, c_k\}$ is dominating set of $N_k$, where $|X| = \frac{n}{4}$. This, together with the lower bound $\gamma(N_k)\geq \frac{n}{4}$, imply that $X$ is a minimum dominating set for $N_k$, and so, $\gamma(N_k)= \frac{n}{4}$. This result can be combined with Lemma~\ref{lem:order-bound}, implying the following formula for the zero forcing number of diamond-necklaces in terms of the domination number. 
\begin{lem}\label{lem:diamond-formula}
If $G \in \mathcal{N}_{cubic}^*$, then $Z(G) = \gamma(G) + 2$. 
\end{lem}

The following result improves on the statement of Conjecture~\ref{conj:confirmed} when the graph $G$ has no diamond-unit in its $\Delta$-D-partition. 
\begin{lem}\label{lem:2-factor}
If $G$ is a connected and cubic graph with a spanning 2-factor consisting only of triangles, then $Z(G) \leq \gamma(G) + 1$.
\end{lem}
\proof Let $G$ be a connected and cubic graph with a spanning 2-factor consisting only of triangles. Thus, $G$ is connected, cubic, and claw-free, where every unit in the $\Delta$-D-partition of $V(G)$ is a triangle-unit. Thus, $G$ is diamond-free. We next define the \emph{contraction multigraph} of $G$, denoted $M_G$, to be the multigraph whose vertices correspond to the triangle-units in the $\Delta$-D-partition of $G$, and where two vertices of $M_G$ are joined by the number of edges joining the corresponding triangle-units in $G$. By construction, $M_G$ has no loops but does possibly contain multiedges. Note that $n(M_G)$ is precisely the number of triangle-units in $G$. Since $G$ contains at least two triangle-units, it must be the case that $n(M_G) \geq 2$. Moreover, each vertex of $M_G$ has degree three, and so, $M_G$ is a cubic multigraph. For each $v \in V(M_G)$, let $V(T_v) = \{a_v, b_v, c_v\}$ denote the vertices of the triangle-unit in $G$ associated with $v$. Further, we say that $T_v$ is derived from $v$ in $M_G$. 

Let $X \subseteq V(G)$ be a minimum (independent) dominating set of $G$ and let $Y = V(G)\setminus X$. Thus, every vertex of $G$ not in $X$ has a neighbor in $X$, and no two vertices in $X$ are endpoints of the same edge in $G$; also, $|X| = \gamma(G)$. Further note that if $T$ is a triangle-unit in the $\Delta$-D-partition of $V(G)$, then either $T$ contains precisely one vertex from $X$ (recall $X$ is an independent set), or each of the triangle-units adjacent with $T_v$ contains a vertex from $X$. Since $G$ is a cubic graph and every vertex not in $X$ has a neighbor in $X$, next observe $\Delta(G[Y]) = 2$. Therefore, $G[Y]$ is either a path, a cycle, or a disjoint union of paths and cycles. In the following two initialization procedures, we define how to color a set of blue vertices using the dominating set $X$, from which we construct a zero forcing set of $G$. 

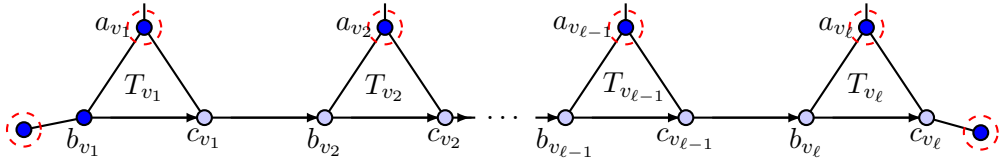
\begin{figure}[htb]
 \begin{center}
\begin{tikzpicture}[scale=.8,style=thick,x=1cm,y=1cm]
\def\vr{3.5pt}

\path (-1.5, 0.0) coordinate (d11);
\path (-2.5, -0.2) coordinate (d111);
\path (-0.5, 1.5) coordinate (d12);
\path (-0.5, 1.9) coordinate (d122);
\path (.5, 0.0) coordinate (d13);
\path (.25, -0.2) coordinate (d133);

\path (2.5, 0.0) coordinate (d21);
\path (2.25, -0.25) coordinate (d211);
\path (3.5, 1.5) coordinate (d22);
\path (3.5, 1.9) coordinate (d222);
\path (4.5, 0.0) coordinate (d23);
\path (5.0, 0.0) coordinate (d233);

\path (6.0, 0.0) coordinate (d311);
\path (6.5, 0.0) coordinate (d31);
\path (7.5, 1.5) coordinate (d32);
\path (7.5, 1.9) coordinate (d322);
\path (8.5, 0.0) coordinate (d33);
\path (8.75, -0.25) coordinate (d333);

\path (10.5, 0.0) coordinate (d41);
\path (10.5, -0.25) coordinate (d411);
\path (11.5, 1.5) coordinate (d42);
\path (11.5, 1.9) coordinate (d422);
\path (12.5,0.0) coordinate (d43);
\path (13.4, -0.25) coordinate (d433);


\draw (d11) -- (d12);
\draw (d11) -- (d111);
\draw (d12) -- (d13);
\draw (d12) -- (d122);
\draw (d11) -- (d13);
\draw (d13) -- (d21);

\draw (d21) -- (d22);
\draw (d22) -- (d23);
\draw (d22) -- (d222);
\draw (d21) -- (d23);
\draw (d23) -- (d233);

\draw (d31) -- (d311);
\draw (d31) -- (d32);
\draw (d32) -- (d33);
\draw (d32) -- (d322);
\draw (d31) -- (d33);
\draw (d33) -- (d41);

\draw (d41) -- (d42);
\draw (d42) -- (d43);
\draw (d42) -- (d422);
\draw (d41) -- (d43);
\draw (d43) -- (d433);

\draw[-latex, shorten >=.65mm, shorten <=3.5mm, line width=0.25mm] (d11) -- (d13);
\draw[-latex, shorten >=.65mm, shorten <=3.5mm, line width=0.25mm] (d13) -- (d21);
\draw[-latex, shorten >=.65mm, shorten <=3.5mm, line width=0.25mm] (d21) -- (d23);
\draw[-latex, shorten >=.65mm, shorten <=3.5mm, line width=0.25mm] (d23) -- (d233);
\draw[-latex, shorten >=.65mm, shorten <=3.5mm, line width=0.25mm] (d311) -- (d31);
\draw[-latex, shorten >=.65mm, shorten <=3.5mm, line width=0.25mm] (d31) -- (d33);

\draw[-latex, shorten >=.65mm, shorten <=3.5mm, line width=0.25mm] (d33) -- (d41);
\draw[-latex, shorten >=.65mm, shorten <=3.5mm, line width=0.25mm] (d41) -- (d43);

\draw (d11) [fill=blue] circle (\vr);
\draw (d111) [fill=blue] circle (\vr);
\draw[dashed, red] (d111) circle (8pt);
\draw (d12) [fill=blue] circle (\vr);
\draw[dashed, red] (d12) circle (8pt);
\draw (d13) [fill=blue!20] circle (\vr);

\draw (d21) [fill=blue!20] circle (\vr);
\draw (d22) [fill=blue] circle (\vr);
\draw[dashed, red] (d22) circle (8pt);
\draw (d23) [fill=blue!20] circle (\vr);

\draw (d31) [fill=blue!20] circle (\vr);
\draw (d32) [fill=blue] circle (\vr);
\draw[dashed, red] (d32) circle (8pt);
\draw (d33) [fill=blue!20] circle (\vr);

\draw (d41) [fill=blue!20] circle (\vr);
\draw (d433) [fill=blue] circle (\vr);
\draw[dashed, red] (d433) circle (8pt);
\draw (d42) [fill=blue] circle (\vr);
\draw[dashed, red] (d42) circle (8pt);
\draw (d43) [fill=blue!20] circle (\vr);

\node at (-.5, 0.5) {$T_{v_1}$};
\node at (3.5, 0.5) {$T_{v_2}$};
\node at (7.65, 0.5) {$T_{v_{\ell-1}}$};
\node at (11.5, 0.5) {$T_{v_{\ell}}$};

\node at (d11) [below] {$b_{v_1}$};
\node at (d12) [left] {$a_{v_1}$};
\node at (d13) [below] {$c_{v_1}$};

\node at (d21) [below] {$b_{v_2}$};
\node at (d22) [left] {$a_{v_2}$};
\node at (d23) [below] {$c_{v_2}$};

\node at (d31) [below] {$b_{v_{\ell-1}}$};
\node at (d32) [left] {$a_{v_{\ell-1}}$};
\node at (d33) [below] {$c_{v_{\ell-1}}$};

\node at (d41) [below] {$b_{v_{\ell}}$};
\node at (d42) [left] {$a_{v_{\ell}}$};
\node at (d43) [below] {$c_{v_{\ell}}$};

\node at (5.5, 0.0) {$\dots$};

\end{tikzpicture}
\end{center}
\caption{The the initialization procedure when $G[Y]$ contains a path component. The initial set $S$ of blue colored vertices in $G$ shown with dark blue color, whereas forcing steps and color changes indicated by directed arrows and light blue colored vertices; dominating vertices indicated by red dashed circles.}
\label{fig:initialization-path}
\end{figure}

\noindent\textbf{Path Initialization.} If $G[Y]$ contains at least one path component, let 
\[
H_G: b_{v_1}a_{v_2}\dots b_{v_{\ell}}c_{v_{\ell}}
\]
be a maximum path component of $G[Y]$, and note that $H_G$ corresponds to the 
path 
\[
H_M : v_1, v_2, \dots, v_{\ell},
\]
in the multigraph $M_G$. Next, color blue the vertices in $S = X \cup \{b_{v_1}\}$ and color white the vertices in $V(G)\setminus S$. Under this coloring of the vertices in $V(G)$, observe that the only white-colored neighbor of $b_{v_1}$ is $c_{v_1}$, and so, the zero forcing process may begin by $b_{v_1}$ forcing $c_{v_1}$ to become colored blue. After $c_{v_1}$ becomes colored blue, then the only possibly white-colored neighbor of $c_{v_1}$ would be $b_{v_{2}}$, and so, $c_{v_1}$ may then $b_{v_{2}}$ to become colored blue. This process continues until $b_{v_{\ell}}$ become colored blue, and thereafter forces $c_{v_{\ell}}$ to become colored blue; see Figure~\ref{fig:initialization-path}. Thus, all vertices in the path component $H_G$ become colored blue, and more so, the derived triangle-units from the vertices of $H_M$ contain only blue-colored vertices. If $G[Y]$ does not contain a path component, move to the next initialization procedure.

\noindent\textbf{Cycle Initialization.} If $G[Y]$ does not contain a path component, then let 
\[
H_G: b_{v_1}a_{v_2}\dots b_{v_{\ell}}c_{v_{\ell}b_{v_1}},
\]
be a largest cycle component of $G[Y]$, where 
\[
H_M: v_1v_2 \dots v_{\ell}v_1,
\]
is a corresponding cycle (allowing 2-cycles) in the multigraph $M_G$. Since $G$ is cubic, claw-free, and diamond-free, a vertex from the dominating set $X$ of $G$ will dominate at most two vertices on the cycle component $H_G$, and if so, then these two vertices must be neighbors on the component $H_G$. Moreover, if a vertex of $w \in X$ dominates exactly one vertex in $H_G$, say $b \in V(H)$, then $H$ is the cycle $C_3$, for otherwise, $v$ and its two neighbors in $H$ together with $w$ would induce the claw $K_{1, 3}$, a contradiction. Thus, since $n(H) \geq 3$, at least one pair of vertices in $H_G$, say $v_1$ and $v_{\ell}$ (renaming if need be), which are dominated by distinct vertices in $X$. Thus, $a_{v_1}$ is the unique vertex dominating vertices in the triangle-unit derived from $v_1$, and $a_{v_{\ell}}$ is the unique vertex dominating the triangle-unit derived from $v_{\ell}$. 
Next color blue the vertices in $S = (X \setminus \{ a_{v_{\ell}} \}) \cup \{b_{v_1}, c_{v_1} \}$ and color white the vertices in $V(G)\setminus S$. 

\begin{figure}[htb]
\begin{center}
\begin{tikzpicture}[scale=.8,style=thick,x=1cm,y=1cm]
\def\vr{3.5pt}

\path (-1.5, 0.0) coordinate (d11);
\path (-2, -0.2) coordinate (d111);
\path (-0.5, 1.5) coordinate (d12);
\path (-0.5, 1.9) coordinate (d122);
\path (.5, 0.0) coordinate (d13);
\path (.25, -0.2) coordinate (d133);

\path (2.5, 0.0) coordinate (d21);
\path (2.25, -0.25) coordinate (d211);
\path (3.5, 1.5) coordinate (d22);
\path (3.5, 1.9) coordinate (d222);
\path (4.5, 0.0) coordinate (d23);
\path (5.0, 0.0) coordinate (d233);

\path (6.0, 0.0) coordinate (d311);
\path (6.5, 0.0) coordinate (d31);
\path (7.5, 1.5) coordinate (d32);
\path (7.5, 1.9) coordinate (d322);
\path (8.5, 0.0) coordinate (d33);
\path (8.75, -0.25) coordinate (d333);

\path (10.5, 0.0) coordinate (d41);
\path (10.25, -0.25) coordinate (d411);
\path (11.5, 1.5) coordinate (d42);
\path (11.5, 1.9) coordinate (d422);
\path (12.5,0.0) coordinate (d43);
\path (12.9, -0.25) coordinate (d433);


\draw (d11) -- (d12);
\draw (d12) -- (d13);
\draw (d12) -- (d122);
\draw[->, bend right, >=latex] (d11) to (d43);
\draw (d11) -- (d13);
\draw[-latex, shorten >=.65mm, shorten <=3.5mm, line width=0.25mm] (d12) -- (d13);
\draw (d13) -- (d21);
\draw[-latex, shorten >=.65mm, shorten <=3.5mm, line width=0.25mm] (d13) -- (d21);

\draw (d21) -- (d22);
\draw (d22) -- (d23);
\draw[-latex, shorten >=.65mm, shorten <=3.5mm, line width=0.25mm] (d21) -- (d23);
\draw (d22) -- (d222);
\draw (d21) -- (d23);
\draw (d23) -- (d233);
\draw[-latex, shorten >=0.25mm, shorten <=3.5mm, line width=0.25mm] (d23) -- (d233);

\draw (d31) -- (d311);
\draw[-latex, shorten >=0.25mm, shorten <=3.5mm, line width=0.25mm] (d311) -- (d31);
\draw (d31) -- (d32);
\draw (d32) -- (d33);
\draw (d32) -- (d322);
\draw (d31) -- (d33);
\draw[-latex, shorten >=.65mm, shorten <=3.5mm, line width=0.25mm] (d31) -- (d33);
\draw (d33) -- (d41);
\draw[-latex, shorten >=.65mm, shorten <=3.5mm, line width=0.25mm] (d33) -- (d41);

\draw (d41) -- (d42);
\draw (d42) -- (d43);
\draw (d42) -- (d422);
\draw (d41) -- (d43);
\draw[-latex, shorten >=.65mm, shorten <=3.5mm, line width=0.25mm] (d41) -- (d42);

\draw (d11) [fill=blue] circle (\vr);
\draw (d12) [fill=blue] circle (\vr);
\draw[dashed, red] (d12) circle (8pt);
\draw (d13) [fill=blue] circle (\vr);

\draw (d21) [fill=blue!20] circle (\vr);
\draw (d22) [fill=blue] circle (\vr);
\draw[dashed, red] (d22) circle (8pt);
\draw (d23) [fill=blue!20] circle (\vr);

\draw (d31) [fill=blue!20] circle (\vr);
\draw (d32) [fill=blue] circle (\vr);
\draw[dashed, red] (d32) circle (8pt);
\draw (d33) [fill=blue!20] circle (\vr);

\draw (d41) [fill=blue!20] circle (\vr);
\draw (d42) [fill=blue!20] circle (\vr);
\draw[dashed, red] (d42) circle (8pt);
\draw (d43) [fill=blue!20] circle (\vr);

\node at (-.5, 0.5) {$T_{v_1}$};
\node at (3.5, 0.5) {$T_{v_2}$};
\node at (7.65, 0.5) {$T_{v_{\ell-1}}$};
\node at (11.5, 0.5) {$T_{v_{\ell}}$};

\node at (d11) [below] {$b_{v_1}$};
\node at (d12) [left] {$a_{v_1}$};
\node at (d13) [below] {$c_{v_1}$};

\node at (d21) [below] {$b_{v_2}$};
\node at (d22) [left] {$a_{v_2}$};
\node at (d23) [below] {$c_{v_2}$};

\node at (d31) [below] {$b_{v_{\ell-1}}$};
\node at (d32) [left] {$a_{v_{\ell-1}}$};
\node at (d33) [below] {$c_{v_{\ell-1}}$};

\node at (d41) [below] {$b_{v_{\ell}}$};
\node at (d42) [left] {$a_{v_{\ell}}$};
\node at (d43) [below] {$c_{v_{\ell}}$};

\node at (5.5, 0.0) {$\dots$};

\end{tikzpicture}
\end{center}
\caption{The the initialization procedure when $G[Y]$ does not contain a path component. The initial set $S$ of blue colored vertices in $G$ shown with dark blue color, whereas forcing steps and color changes indicated by directed arrows and light blue colored vertices; dominating vertices indicated by red dashed circles.}
\label{fig:intialization-cycle}
\end{figure}
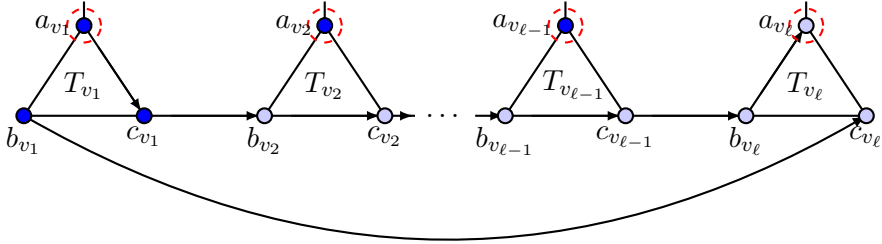

Under this coloring of vertices in $V(G)$, observe that the only white-colored neighbor of $b_{v_1}$ is $c_{v_{\ell}}$ and the only white colored neighbor of $c_{v_1}$ is $b_{v_2}$. Thus, the zero forcing process may begin by having $b_{v_1}$ force $c_{v_{\ell}}$ to become colored blue, and also having $c_{v_1}$ force $b_{v_2}$ to become colored blue. Thereafter, the only white-colored neighbor of $c_{v_1}$ would be $b_{v_2}$, and so, $c_{v_1}$ may then force $b_{v_2}$ to become colored blue, and this process continues until $b_{v_{\ell}}$ becomes colored blue. Then after $b_{v_{\ell}}$ becomes colored blue, recall that $b_{v_1}$ had previously forced $c_{v_{\ell}}$ to become colored blue. Thus, the only currently white-colored neighbor of $b_{v_{\ell}}$ is $a_{v_{\ell}}$, and so, $b_{v_{\ell}}$ would now force $a_{v_{\ell}}$ to become colored blue; see Figure~\ref{fig:intialization-cycle} for an illustration. Thus, all vertices in the cycle $H_G$ become colored blue, and more so, the derived triangle-units from the vertices of $H_M$ contain only blue-colored vertices. 

Let $S \subseteq V(G)$ be the set of vertices obtained by either of the above initialization techniques. and thereafter, allow the zero forcing process to start and continue until no further color changes are possible. By the construction of the initial set of blue colored vertices in $S$, observe the following properties:
\begin{itemize}
    \item $|S| = |X| + 1$.
    \item All vertices in the dominating set $X$ of $G$ are currently blue-colored. 
    \item Each of triangle-units derived from the subgraph $H_M$ contain blue-colored vertices in $G$. 
    \item $X$ is reachable from $S$.
\end{itemize}

If all vertices of $G$ have become blue-colored, then $S$ is a zero forcing set of $G$, and so, 
\[
Z(G) \leq |S| = |X| + 1 = \gamma(G) + 1.
\]
Hence, we may assume that not all vertices of $G$ are reachable by $S$ since the desired result follows otherwise.  

\begin{figure}[htb]
\begin{center}
\begin{tikzpicture}[scale=.8,style=thick,x=1cm,y=1cm]
\def\vr{3.5pt}

\path (-1.5, 0.0) coordinate (d11);
\path (-1.75, -0.25) coordinate (d111);
\path (-0.5, 1.5) coordinate (d12);
\path (-0.5, 1.9) coordinate (d122);
\path (.5, 0.0) coordinate (d13);
\path (.75, -0.25) coordinate (d133);

\path (2.5, 0.0) coordinate (d21);
\path (2.25, -0.25) coordinate (d211);
\path (3.5, 1.5) coordinate (d22);
\path (3.5, 1.9) coordinate (d222);
\path (4.5, 0.0) coordinate (d23);
\path (4.75, -0.25) coordinate (d233);

\path (6.5, 0.0) coordinate (d31);
\path (6.25, -0.25) coordinate (d311);
\path (7.5, 1.5) coordinate (d32);
\path (7.5, 1.9) coordinate (d322);
\path (8.5, 0.0) coordinate (d33);
\path (8.75, -0.25) coordinate (d333);

\path (10.5, 0.0) coordinate (d41);
\path (10.25, -0.25) coordinate (d411);
\path (11.5, 1.5) coordinate (d42);
\path (11.5, 1.9) coordinate (d422);
\path (12.5,0.0) coordinate (d43);
\path (12.75, -0.25) coordinate (d433);

\filldraw[fill=blue!20] (d41) -- (d42) -- (d43) -- cycle;

\draw (d11) -- (d12);
\draw (d11) -- (d111);
\draw (d12) -- (d13);
\draw (d12) -- (d122);
\draw (d11) -- (d13);
\draw (d13) -- (d133);

\draw (d21) -- (d22);
\draw (d21) -- (d211);
\draw (d22) -- (d23);
\draw (d22) -- (d222);
\draw (d21) -- (d23);
\draw (d23) -- (d233);

\draw (d31) -- (d32);
\draw (d31) -- (d311);
\draw (d32) -- (d33);
\draw (d32) -- (d322);
\draw (d31) -- (d33);
\draw (d33) -- (d333);

\draw (d41) -- (d42);
\draw (d41) -- (d411);
\draw (d42) -- (d43);
\draw (d42) -- (d422);
\draw (d41) -- (d43);
\draw (d43) -- (d433);

\draw (d11) [fill=white] circle (\vr);
\draw (d12) [fill=white] circle (\vr);
\draw (d13) [fill=white] circle (\vr);

\draw (d21) [fill=white] circle (\vr);
\draw (d22) [fill=blue] circle (\vr);
\draw (d23) [fill=white] circle (\vr);

\draw (d31) [fill=white] circle (\vr);
\draw (d32) [fill=blue] circle (\vr);
\draw (d33) [fill=blue] circle (\vr);

\draw (d41) [fill=blue] circle (\vr);
\draw (d42) [fill=blue] circle (\vr);
\draw (d43) [fill=blue] circle (\vr);


\node at (-.5, -1.0) {$\frac{0}{3}$-blue colored};
\node at (3.5, -1.0) {$\frac{1}{3}$-blue colored};
\node at (7.5, -1.0) {$\frac{2}{3}$-blue colored};
\node at (11.5, -1.0) {$\frac{3}{3}$-blue colored};

\end{tikzpicture}
\end{center}
\caption{The different types of $\frac{k}{3}$-blue colorings of the triangle-units in $G$ as given in the proof of Lemma~\ref{lem:2-factor}.}
\label{fig:k/3-coloring}
\end{figure}
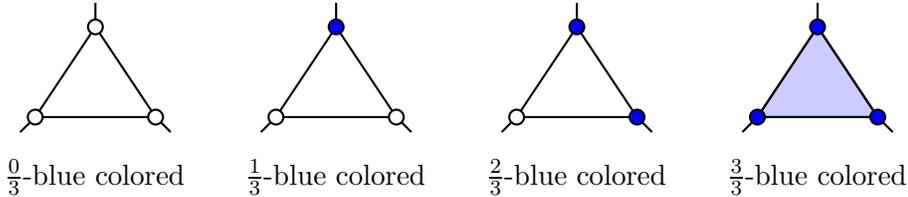

Suppose the zero forcing process stopped before all vertices in $V(G)$ become blue-colored. By constructing $S$, we know that all vertices in $P_G$ are currently blue-colored, and further, that any derived triangle-unit from a vertex in $H_M$ will only contain blue-colored vertices in $G$. We will now track color changes in $G$ through associated color changes in the triangle-units of $G$, which we now define in terms of the multigraph $M_G$. For a given blue and white coloring of $V(G)$ and a vertex $v \in V(M_G)$, we say that $v$ (and also $T_v$) is \emph{$\frac{k}{3}$-blue-colored} whenever $V(T_v)$ contains exactly $k$ blue currently colored vertices in $G$, where $k \in [3]$; see Figure~\ref{fig:k/3-coloring}. Thus, there is a duality between blue and white colorings of $V(G)$ and $\frac{k}{3}$-blue colorings of the triangle-units in $G$. That is, a coloring of $V(G)$ implies a coloring of the triangle-units of $G$ and vice versa. We define a \emph{blue-white-open} edge to be an edge of $G$ with a blue-colored endpoint in a $\frac{3}{3}$-blue-colored triangle-unit and a white-colored endpoint in an adjacent triangle-unit. 

The following two claims significantly simplify the monitoring of color changes in $G$ regarding the derived triangle units from vertices in $M_G$. 
\begin{claim}\label{claim:1-blue-one}
Let $v$ and $w$ be two adjacent vertices of $M_G$ adjacent by at least one blue-white-open edge in $G$. If $v$ is $\frac{3}{3}$-blue-colored and $w$ is $\frac{1}{3}$-blue-colored, then $w$ becomes $\frac{3}{3}$-blue-colored by zero forcing color changes in $G$. 
\end{claim}
\proof Let $v$ and $w$ be two adjacent vertices of $M_G$ which are adjacent by at least one blue-white-open edge in $G$. Recalling our choice of notation, let $T_v = \{a_v, b_v, c_v\}$ and $V(T_w) = \{a_w, b_w, c_w\}$ denote the vertices of the derived triangle-units from $v$, and $w$, respectively. Suppose $v$ is $\frac{3}{3}$-blue colored and $w$ is $\frac{1}{3}$-blue colored, where $a_v a_w$ is a blue-white-open edge connecting $v$ to $w$. Since $w$ is $\frac{1}{3}$-blue colored, this supposition implies exactly one of the vertices in $\{b_w, c_w\}$ is currently blue-colored, say $c_w$. Moreover, since $v$ is $\frac{3}{3}$-blue-colored, $a_w$ is the only white-colored neighbor of the blue-colored vertex $a_v$ in $G$. Thus, $a_v$ may force $a_w$ to become blue-colored. After $a_w$ becomes blue-colored, the only white-colored neighbor of $a_w$ would be $b_w$, and so, $a_w$ may then force $b_w$ to become blue-colored. Therefore, by these aforementioned color changes in $G$, the vertex $w$ in $M_G$ becomes $\frac{3}{3}$-blue colored. \smallqed 

\begin{claim}\label{claim:1-blue-two}
Let $v$ and $w$ be two adjacent vertices of $M_G$. If $v$ is $\frac{3}{3}$-blue-colored and $w$ is $\frac{2}{3}$-blue-colored, then $w$ will become $\frac{3}{3}$-blue-colored by way of zero forcing color changes in $G$. 
\end{claim}
\proof Let $v$ and $w$ be two adjacent vertices of $M_G$, and suppose $v$ is $\frac{3}{3}$-blue-colored and $w$ is $\frac{2}{3}$-blue-colored. Recalling our choice of notation, let $T_v = \{a_v, b_v, c_v\}$ and $V(T_w) = \{a_w, b_w, c_w\}$ denote the vertices of the derived triangle-units from $v$, and $w$, respectively. Let $a_va_w$ be an edge of $G$ connecting the triangle-units $T_v$ and $T_w$. Since $v$ is $\frac{3}{3}$-blue-colored, every vertex in $V(T_v)$ is blue-colored in $G$. Since $w$ is $\frac{2}{3}$-blue-colored, exactly two of the vertices in $V(T_w)$ are blue-colored in $G$. If $a_w$ is white-colored, then $b_w$ and $c_w$ are blue-colored in $G$. With this supposition, $a_w$ would then be the only white-colored neighbor of $a_v$, and so, $a_v$ would then force $a_w$ to become blue-colored, resulting in $w$ becoming $\frac{3}{3}$-blue-colored. If $a_w$ were currently blue-colored, then exactly one of the vertices in $\{b_w, c_w\}$ is blue-colored, say $c_w$. Thus, since $a_v$ is blue-colored, the only white-colored neighbor of $a_w$ would be $b_w$, and so, $a_w$ would then force $b_w$ to become blue-colored. Hence, $w$ becomes $\frac{3}{3}$-blue-colored. \smallqed 

With the above claims proven, we now return to the current state of blue and white colored vertices in $V(G)$. Since each vertex in $H_M$ is $\frac{3}{3}$-blue-colored, we know the current set of $\frac{3}{3}$-blue-colored vertices is nonempty. Let $B$ be the current set of $\frac{3}{3}$-blue-colored vertices in $M_G$. Since $M_G$ is a connected multigraph, and not all vertices of $G$ are blue-colored, at least one vertex of $M_G$ must exist that is (noninclusive) adjacent to a vertex in $B$. Let $w \in V(M_G)$ be a vertex adjacent to the set $B$ in $M_G$. Note that $w$ cannot be $\frac{0}{3}$-blue-colored, since if so, one of the vertices in $B$ would have forced onto $T_w$ making $w$ at least $\frac{1}{3}$-blue-colored. If $w$ were $\frac{2}{3}$-blue-colored, then by Claim~\ref{claim:1-blue-two}, $w$ would become $\frac{3}{3}$-blue-colored, and hence belong to $B$. Thus, $w$ must be $\frac{1}{3}$-blue-colored. However, by Claim~\ref{claim:1-blue-one}, if $w$ were adjacent with $B$ by at least one blue-white-open edge in $G$, then $w$ would become $\frac{3}{3}$-blue-colored, and hence belong to $B$. Hence, $w$ is $\frac{1}{3}$-blue-colored, and no edge connecting $w$ to $B$ is blue-white-open. Both endpoints of exactly one edge connecting $w$ to $B$ are blue-colored. 

Since we are assuming that $S$ is not a zero forcing set $G$, we now define a collection of \emph{extension rules} that modify the set $S$ so that after each application of the said rule, more vertices become colored blue; also we maintain the cardinality of our original set $S$. For each rule we introduce, we will formally state it and then prove a claim for its validity. For example, see the first simple rule given below. 

\begin{swap}\label{swap:easy}
Let $w \notin B$ be adjacent with $v \in B$ in $M_G$. If $a_va_w$ is the edge connecting $v$ to $w$ in $G$ and $a_v$ did not force during the zero forcing process previously, then
\[
S \leftarrow (S \setminus \{a_w\}) \cup \{b_w\}.
\]
\end{swap}

\begin{claim}\label{claim:easy}
Extension Rule~\ref{swap:easy} is a valid-$z$-proper-extension rule on $S$.
\end{claim}

\proof Let $w \notin B$ is adjacent with $v \in B$ in $M_G$. Recalling our choice of notation, let $T_v = \{a_v, b_v, c_v\}$ and $V(T_w) = \{a_w, b_w, c_w\}$ denote the vertices of the derived triangle-units from $v$, and $w$, respectively. Next, suppose $a_va_w$ is the edge connecting $v$ to $w$ in $G$, and, further, that $a_v$ did not force during the zero forcing process previously. Note $a_w$ is blue-colored. Let $S' = (S \setminus \{a_w\}) \cup \{b_w\}$, and so, $|S'| = |S|$. Thus, it remains to show $N_G^z[S] \subset N_G^z[S']$. Since $a_v$ did not force a color change during the zero forcing process, and since $b_w$ and $c_w$ are white-colored, removing $a_w$ does not affect any color changes that occurred by applying the zero forcing process on $S$ in $G$. Thus, all color changes due to $S$ will also occur with $S'$. In particular, at some point in the zero forcing process, each vertex in the triangle-unit $T_v$ is colored blue, after which $a_w$ would be the only white-colored neighbor of $a_v$. Therefore, $a_v$ may then force $a_w$ to become colored blue. Thus, $N_G^z[S] \subseteq N_G^z[S']$. However, since $b_w \notin  N_G^z[S]$ and $b_w \in  N_G^z[S']$, we establish $N_G^z[S] \subset N_G^z[S']$, and the desired result follows. \smallqed 

By Claim~\ref{claim:easy}, we may apply Extension Rule~\ref{swap:easy} as many times as needed, always producing a new initial set of $S$ of blue-colored vertices in $G$ and always maintaining $|S| = |X| + 1$. If we eventually arrive at a set $S$ of blue-colored vertices so that all of $V(G)$ becomes blue-colored, then $S$ is a zero forcing set of $G$. Thus, 
\[
Z(G) \leq |S| = |X| + 1 = \gamma(G) + 1. 
\]
Therefore, we will assume that not all of $V(G)$ has become blue-colored since otherwise, the desired result will follow. 

\begin{figure}[htb]
 \begin{center}
\begin{tikzpicture}[scale=.8,style=thick,x=1cm,y=1cm]
\def\vr{3.5pt}

\path (1.0, 3.0) coordinate (av);
\path (1.0, 1.0) coordinate (aw);
\path (-0.5, -1.0) coordinate (bw);
\path (2.5, -1.0) coordinate (cw);

\draw (aw) -- (bw);
\draw (aw) -- (cw);
\draw (bw) -- (cw);
\draw (av) -- (aw);
\draw[-latex, shorten >=.65mm, shorten <=3.5mm, line width=0.25mm] (av) -- (aw);

\draw (av) [fill=blue] circle (\vr);
\draw (aw) [fill=blue!20] circle (\vr);
\draw[dashed, red] (av) circle (25pt);
\draw (bw) [fill=white] circle (\vr);
\draw (cw) [fill=white] circle (\vr);

\node at (av) [above] {$a_v$};
\node at (aw) [left] {$a_z$};
\node at (bw) [below] {$b_z$};
\node at (cw) [below] {$c_z$};

\node at (1.0, -0.25) {$T_{z}$};
\node at (1, -1.75) {(a)};

\begin{scope}[shift={(6.5,0)}]
    \path (1.0, 3.0) coordinate (av);
    \path (1.0, 1.0) coordinate (aw);
    \path (-0.5, -1.0) coordinate (bw);
    \path (2.5, -1.0) coordinate (cw);
    
    \draw (aw) -- (bw);
    \draw (aw) -- (cw);
    \draw (bw) -- (cw);
    
    \draw (av) -- (aw);
    \draw[-latex, shorten >=.65mm, shorten <=3.5mm, line width=0.25mm] (av) -- (aw);
    
    \draw (av) [fill=blue] circle (\vr);
    \draw (aw) [fill=blue!20] circle (\vr);
    \draw[dashed, red] (aw) circle (25pt);
    \draw (bw) [fill=white] circle (\vr);
    \draw (cw) [fill=white] circle (\vr);

    \node at (av) [above] {$a_v$};
    \node at (aw) [left] {$a_z$};
    \node at (bw) [below] {$b_z$};
    \node at (cw) [below] {$c_z$};
    
    \node at (1.0, -0.25) {$T_{z}$};

    \node at (1, -1.75) {(b)};
\end{scope}

\end{tikzpicture}
\end{center}
\caption{(a) The vertex $a_v \in X$ forces $a_z \notin X$ to become blue-colored. (b) The vertex $a_v \notin X$ forces $a_z \in X$ to become blue-colored. In both instances, the dominating vertex of the configuration is indicated by the red dashed circle.}
\label{fig:cases}
\end{figure}
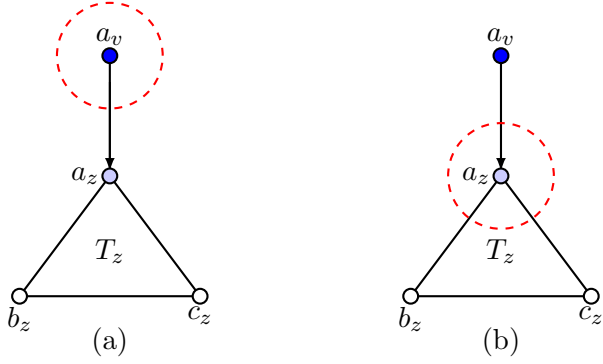

Let $z$ be a $\frac{1}{3}$-blue colored vertex of $M_G$ adjacent with the set $B$ and note that Extension Rule~\ref{swap:easy} implies the blue-colored vertex of $T_z$ was not originally colored in the set $S$. Thus, the blue-colored vertex in $T_z$ must have originally been white-colored and become blue-colored by the zero forcing process in $G$. Let $a_z$ be the vertex of $T_z$ that became blue-colored during the zero forcing process in $G$, and so, $a_z$ is the endpoint of some edge whose other endpoint is in $B$. There are two possible cases for the vertices in the triangle-unit $T_z$ concerning the dominating set $X$ of $G$; either $V(T_z) \cap X \neq \emptyset$ or $V(T_z) \cap X = \emptyset$; see Figure~\ref{fig:cases}. If $a_z \in X$, we are in the configuration shown by Figure~\ref{fig:cases}~(a). In this case, recall that all vertices of $X$ are currently colored blue (either in $S$ or by some forcing step). Thus, $b_z$ and $c_z$ belong to some component of $G[Y]$, say $H_G^z$, and this component is either a path or a cycle. We are considering both possibilities now. 

Suppose first that $a_v \in X$ forces a neighbor $a_z$ to become blue-colored. Since $b_z$ and $c_z$ must be dominated by a vertex in $X$, $b_z$ and $c_z$ have a neighbor in $X$, say $a_x$ and $a_w$, respectively. Since all vertices of $X$ are reachable from $S$, we remark that $a_x$ and $a_w$ are currently blue-colored.

We next introduce an Extension Rule for handling this case, provided Extension Rule~\ref{swap:easy} has been applied as many times as possible prior. 
\begin{swap}\label{swap:empty-triangle}
Let $z \notin B$ be adjacent with $v \in B$ in $M_G$. If $a_v \in X$ forces $a_z \notin X$ to become blue-colored, then
\[
S \leftarrow (S \setminus \{a_w\}) \cup \{b_z\}.
\]
\end{swap}

\begin{claim}\label{claim:empty-triangle}
Extension Rule~\ref{swap:empty-triangle} is a valid-$z$-proper-extension rule on $S$
\end{claim}
\proof Note that since the $\Delta$-D-partition of $V(G)$ contains no diamond-unit, $a_x \neq a_w$. Let $S' = (S \setminus \{a_w\}) \cup \{b_z\}$, and so, $|S'| = |S|$. Thus, it remains to show $N_G^z[S] \subset N_G^z[S']$. If any neighbor of $a_w$ forced during the zero forcing process under the coloring starting with $S$, then all vertices of the triangle-unit $T_y$ would necessarily be colored blue, which implies that $y$ was $\frac{3}{3}$-blue-colored. Thus, $y$ may force onto $T_z$ across the blue-white-open edge $a_wc_z$. By Claim~\ref{claim:1-blue-one}, $T_{z}$ would then become $\frac{3}{3}$-blue-colored, a contradiction since $b_z$ and $c_z$ are currently white-colored. Thus, no neighbor of $a_w$ forces a color change during the zero forcing process starting with $S$. Therefore, all color changes previously due to $S$ will also occur by starting the zero forcing process with $S'$. In particular, at some point of the zero forcing process, $a_v$ will force $a_z$ to become blue-colored. Thereafter, the only white-colored neighbor of $a_z$ would be $c_z$, so $a_z$ may force $c_z$ to become blue-colored. After $c_z$ becomes blue-colored, the only white-colored neighbor of $c_z$ would be $a_w$, and so, $c_w$ may force $a_w$ to become blue-colored. Thus, $N_G^z[S'] = N_G^z[S] \cup \{b_z, c_z\}$. Hence, $N_G^z[S] \subset N_G^z[S']$, and the proof of the claim is finished. \smallqed 

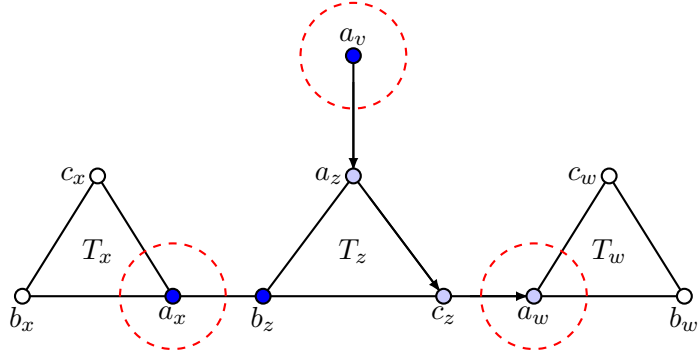
\begin{figure}[htb]
 \begin{center}
\begin{tikzpicture}[scale=.8,style=thick,x=1cm,y=1cm]
\def\vr{3.5pt}

\path (1.0, 3.0) coordinate (av);
\path (1.0, 1.0) coordinate (aw);
\path (-0.5, -1.0) coordinate (bw);
\path (2.5, -1.0) coordinate (cw);

\path (-2.0, -1.0) coordinate (ax);
\path (-4.5, -1.0) coordinate (bx);
\path (-3.25, 1.0) coordinate (cx);

\path (4.0, -1.0) coordinate (ay);
\path (6.5, -1.0) coordinate (by);
\path (5.25, 1.0) coordinate (cy);

\draw (aw) -- (bw);
\draw (aw) -- (cw);
\draw (bw) -- (cw);
\draw (av) -- (aw);
\draw[-latex, shorten >=.65mm, shorten <=3.5mm, line width=0.25mm] (av) -- (aw);

\draw (bw) -- (ax);
\draw (cw) -- (ay);

\draw (ax) -- (bx);
\draw (ax) -- (cx);
\draw (bx) -- (cx);

\draw (ay) -- (by);
\draw (ay) -- (cy);
\draw (by) -- (cy);

\draw[-latex, shorten >=.65mm, shorten <=3.5mm, line width=0.25mm] (aw) -- (cw);
\draw[-latex, shorten >=.65mm, shorten <=3.5mm, line width=0.25mm] (cw) -- (ay);

\draw (av) [fill=blue] circle (\vr);
\draw (aw) [fill=blue!20] circle (\vr);
\draw[dashed, red] (av) circle (25pt);
\draw (bw) [fill=blue] circle (\vr);
\draw (cw) [fill=blue!20] circle (\vr);

\draw (ax) [fill=blue] circle (\vr);
\draw[dashed, red] (ax) circle (25pt);
\draw (bx) [fill=white] circle (\vr);
\draw (cx) [fill=white] circle (\vr);

\draw (ay) [fill=blue!20] circle (\vr);
\draw[dashed, red] (ay) circle (25pt);
\draw (by) [fill=white] circle (\vr);
\draw (cy) [fill=white] circle (\vr);

\node at (av) [above] {$a_v$};
\node at (aw) [left] {$a_z$};
\node at (bw) [below] {$b_z$};
\node at (cw) [below] {$c_z$};

\node at (ay) [below] {$a_w$};
\node at (by) [below] {$b_w$};
\node at (cy) [left] {$c_w$};

\node at (ax) [below] {$a_x$};
\node at (bx) [below] {$b_x$};
\node at (cx) [left] {$c_x$};

\node at (1.0, -0.25) {$T_{z}$};
\node at (5.25, -0.25) {$T_{w}$};
\node at (-3.25, -0.25) {$T_{x}$};

\end{tikzpicture}
\end{center}
\caption{The vertex $a_v \in X$ forces $a_w \notin X$ to become blue-colored}
\label{fig:empty-triangle-1}
\end{figure}

Apply Extension Rule~\ref{swap:empty-triangle}, and note that (as seen in the proof of Claim~\ref{claim:empty-triangle}), that $x$ may become $\frac{3}{3}$-blue-colored by applying Extension Rule~\ref{swap:easy}, while we do not yet have a rule for $w$. Continue the zero forcing process and apply Extension Rule~\ref{swap:easy} whenever possible until no further color changes are possible. If $w$ does not become $\frac{3}{3}$-blue-colored, then we are in the case as seen by Figure~\ref{fig:cases}~(b), and we proceed to the following two extension rules given below; otherwise $y$ becomes $\frac{3}{3}$-blue-colored and we may once again apply Extension Rule~\ref{swap:empty-triangle}. That is, you are not allowed to apply Extension Rule~\ref{swap:empty-triangle} again unless you have first ensured that all triangle-units containing a \emph{``swapped dominating vertex''} used by an extension rule are $\frac{3}{3}$-blue-colored. Since $a_w$ is a dominating vertex, we note that $b_w$ and $c_w$ must belong to some component of $G[Y]$, say $H_G^w$. 

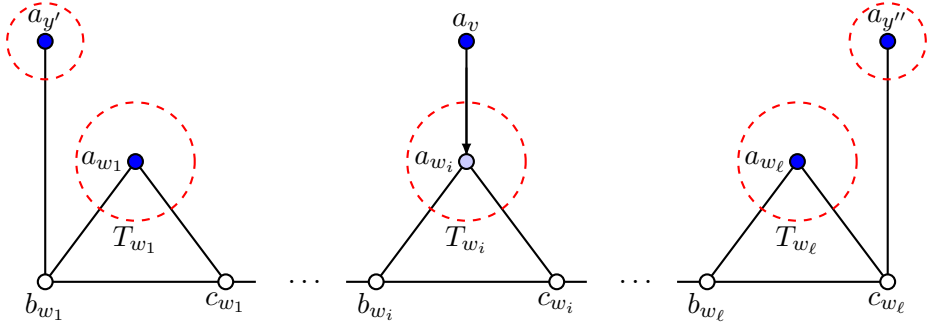
\begin{figure}[htb]
 \begin{center}
\begin{tikzpicture}[scale=.8,style=thick,x=1cm,y=1cm]
\def\vr{3.5pt}

\path (1.0, 1.0) coordinate (aw1);
\path (-0.5, -1.0) coordinate (bw1);
\path (-0.5, 3.0) coordinate (bw11);
\path (2.5, -1.0) coordinate (cw1);
\path (3, -1.0) coordinate (cw11);

\draw (aw1) -- (bw1);
\draw (aw1) -- (cw1);
\draw (bw1) -- (cw1);
\draw (bw1) -- (bw11);
\draw (cw1) -- (cw11);

\draw (aw1) [fill=blue] circle (\vr);
\draw[dashed, red] (aw1) circle (28pt);
\draw (bw1) [fill=white] circle (\vr);
\draw (bw11) [fill=blue] circle (\vr);
\draw[dashed, red] (bw11) circle (18pt);
\draw (cw1) [fill=white] circle (\vr);

\node at (aw1) [left] {$a_{w_1}$};
\node at (bw1) [below] {$b_{w_1}$};
\node at (bw11) [above] {$a_{y'}$};
\node at (cw1) [below] {$c_{w_1}$};

\node at (1.0, -0.3) {$T_{w_1}$};

\node at (3.85, -1) {$\dots$};

\begin{scope}[shift={(5.5,0)}]
    \path (1.0, 3.0) coordinate (av);
    \path (1.0, 1.0) coordinate (awi);
    \path (-0.5, -1.0) coordinate (bwi);
    \path (2.5, -1.0) coordinate (cwi);
    \path (3, -1.0) coordinate (cwi1);

    \path (-1.0, -1.0) coordinate (bwi1);

    \draw (awi) -- (bwi);
     \draw (bwi) -- (bwi1);
    \draw (awi) -- (cwi);
    \draw (bwi) -- (cwi);
    \draw (av) -- (awi);
    \draw (cwi) -- (cwi1);
    \draw[-latex, shorten >=.65mm, shorten <=3.5mm, line width=0.25mm] (av) -- (awi);
    
    \draw (av) [fill=blue] circle (\vr);
    \draw (awi) [fill=blue!20] circle (\vr);
    \draw[dashed, red] (awi) circle (28pt);
    \draw (bwi) [fill=white] circle (\vr);
    \draw (cwi) [fill=white] circle (\vr);
    
    \node at (av) [above] {$a_v$};
    \node at (awi) [left] {$a_{w_i}$};
    \node at (bwi) [below] {$b_{w_i}$};
    \node at (cwi) [below] {$c_{w_i}$};
    
    \node at (1.0, -0.3) {$T_{w_i}$};
    \node at (3.85, -1) {$\dots$};
\end{scope}

\begin{scope}[shift={(11.0,0)}]
    \path (1.0, 1.0) coordinate (awl);
    \path (-0.5, -1.0) coordinate (bwl);
    \path (2.5, -1.0) coordinate (cwl);
    \path (-1.0, -1.0) coordinate (bwl1);
    \path (2.5, 3.0) coordinate (cwl1);

    \draw (awl) -- (bwl);
    \draw (awl) -- (cwl);
    \draw (bwl) -- (cwl);
    \draw (bwl) -- (bwl1);
    \draw (cwl) -- (cwl1);
    
    \draw (awl) [fill=blue] circle (\vr);
    \draw[dashed, red] (awl) circle (28pt);
    \draw (bwl) [fill=white] circle (\vr);
    \draw (cwl) [fill=white] circle (\vr);
    \draw (cwl1) [fill=blue] circle (\vr);
    \draw[dashed, red] (cwl1) circle (18pt);
    
    \node at (awl) [left] {$a_{w_{\ell}}$};
    \node at (bwl) [below] {$b_{w_{\ell}}$};
    \node at (cwl) [below] {$c_{w_{\ell}}$};
    \node at (cwl1) [above] {$a_{y''}$};
    
    \node at (1.0, -0.3) {$T_{w_{\ell}}$};
\end{scope}

\end{tikzpicture}
\end{center}
\caption{A dominating vertex $a_{w_i}$ is forced to become blue-colored and belongs to a path in the multigraph $M_G$.}
\label{fig:path-fix-1}
\end{figure}

Suppose that $H_G^w$ is a path component of $G[Y]$. Denote to the vertices of this component by $H_G^w: b_{w_1}c_{w_1}\dots b_{w_{\ell}}c_{w_{\ell}}$, and denote the associated vertices (triangle-units) in $M_G$ by $H_M^w: w_1w_2\dots w_{\ell}$. Note for some $i \in [\ell]$, $a_{w_i} = a_w$, $b_{w_i} = b_w$, and $c_{w_i} = c_w$. Since $H_G^w$ is a path component of $G[Y]$ and since $G$ is cubic, the enpoints of $H_G^w$ are necessarily adjacent with two vertices from the dominating set $X$ of $G$, say $a_{y'}$ and $a_{w_1}$ for $b_{w_1}$, and $a_{w_{\ell}}$ and $a_{y''}$ for $c_{w_{\ell}}$; see Figure~\ref{fig:path-fix-1} for an illustration of this configuration. The following extension rule addresses this possibility. 

\begin{swap}\label{swap:paths}
Let $a_v$ be a vertex that forced the dominating vertex $a_{w_i}$ in the triangle-unit $T_{w_i}$ with vertex set $V(T_{w_i}) = \{a_{w_i}, b_{w_i}, c_{w_i} \}$. If the component $H_G$ containing $b_{w_i}$ and $c_{w_i}$ is a path, then
\[
S \leftarrow (S \setminus \{a_{y''}\}) \cup \{b_{w_i}\}, 
\]
where $a_{y''}$ is adjacent with $c_{w_{\ell}}$ (and $a_{y''} \notin T_{w_i}$). 
\end{swap}

With the following claim, we establish the validity of Extension Rule~\ref{swap:paths}.
\begin{claim}\label{claim:path-extension}
Extension Rule~\ref{swap:paths} is a valid-$z$-proper-extension rule on $S$.
\end{claim}
\proof Let $S' = (S \setminus \{a_{y''}\}) \cup \{b_{w_i}\}$, and so, $|S'| = |S|$. Thus is remains to show that $N_G^z[S] \subset N_G^z[S']$. If any neighbor of $a_{y''}$ forced during the zero forcing process under the coloring starting with $S$, then all vertices of the triangle-unit $T_{y''}$ would necessarily be blue-colored, which implies that $y''$ was $\frac{3}{3}$-blue-colored. Thus, $y''$ may then forced onto $T_{w_{\ell}}$ across the blue-white-open edge $a_{y''}b_{w_{\ell}}$. By Claim~\ref{claim:1-blue-one}, $T_{w_{\ell}}$ would then become $\frac{3}{3}$-blue-colored. There after, $T_{w_{\ell}}$ could then force onto $T_{w_{\ell - 1}}$ by the blue-white-open edge $b_{w_{\ell}}c_{w_{\ell-1}}$, which in turn would make $w_{\ell - 1}$ a $\frac{3}{3}$-blue-colored vertex. This process would continue until $w_i$ is forced onto, making it a $\frac{3}{3}$-blue-colored vertex, a contradiction since $b_{w_{i}}$ and $c_{w_{i}}$ are currently white-colored. Thus, no neighbor of $a_{y''}$ forces a color change during the zero forcing process starting with $S$. Therefore, all color changes previously due to $S$ will also occur by starting the zero forcing process at $S'$. 

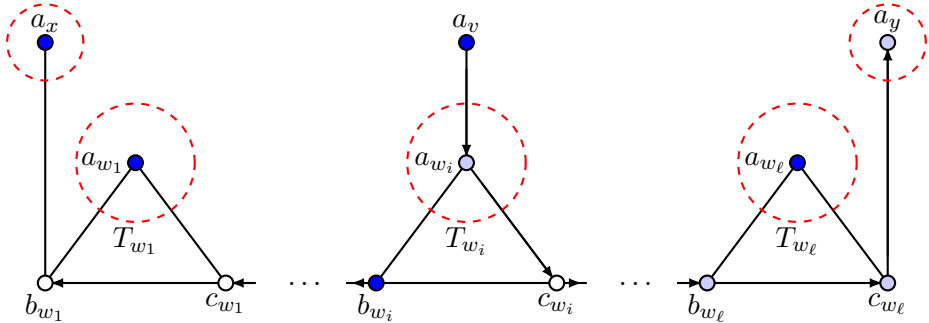
\begin{figure}[htb]
 \begin{center}
\begin{tikzpicture}[scale=.8,style=thick,x=1cm,y=1cm]
\def\vr{3.5pt}

\path (1.0, 1.0) coordinate (aw1);
\path (-0.5, -1.0) coordinate (bw1);
\path (-0.5, 3.0) coordinate (bw11);
\path (2.5, -1.0) coordinate (cw1);
\path (3, -1.0) coordinate (cw11);

\draw (aw1) -- (bw1);
\draw (aw1) -- (cw1);
\draw (bw1) -- (cw1);
\draw (bw1) -- (bw11);
\draw (cw1) -- (cw11);

\draw[-latex, shorten >=.65mm, shorten <=3.5mm, line width=0.25mm] (cw11) -- (cw1);
\draw[-latex, shorten >=.65mm, shorten <=3.5mm, line width=0.25mm] (cw1) -- (bw1);

\draw (aw1) [fill=blue] circle (\vr);
\draw[dashed, red] (aw1) circle (28pt);
\draw (bw1) [fill=white] circle (\vr);
\draw (bw11) [fill=blue] circle (\vr);
\draw[dashed, red] (bw11) circle (18pt);
\draw (cw1) [fill=white] circle (\vr);

\node at (aw1) [left] {$a_{w_1}$};
\node at (bw1) [below] {$b_{w_1}$};
\node at (bw11) [above] {$a_x$};
\node at (cw1) [below] {$c_{w_1}$};

\node at (1.0, -0.3) {$T_{w_1}$};

\node at (3.85, -1) {$\dots$};

\begin{scope}[shift={(5.5,0)}]
    \path (1.0, 3.0) coordinate (av);
    \path (1.0, 1.0) coordinate (awi);
    \path (-0.5, -1.0) coordinate (bwi);
    \path (2.5, -1.0) coordinate (cwi);
    \path (3, -1.0) coordinate (cwi1);

    \path (-1.0, -1.0) coordinate (bwi1);

    \draw (awi) -- (bwi);
     \draw (bwi) -- (bwi1);
    \draw (awi) -- (cwi);
    \draw (bwi) -- (cwi);
    \draw (av) -- (awi);
    \draw (cwi) -- (cwi1);
    \draw[-latex, shorten >=.65mm, shorten <=3.5mm, line width=0.25mm] (av) -- (awi);
    \draw[-latex, shorten >=.65mm, shorten <=3.5mm, line width=0.25mm] (awi) -- (cwi);
    \draw[-latex, shorten >=.65mm, shorten <=3.5mm, line width=0.25mm] (cwi) -- (cwi1);
    \draw[-latex, shorten >=.65mm, shorten <=3.5mm, line width=0.25mm] (bwi) -- (bwi1);
    
    \draw (av) [fill=blue] circle (\vr);
    \draw (awi) [fill=blue!20] circle (\vr);
    \draw[dashed, red] (awi) circle (28pt);
    \draw (bwi) [fill=blue] circle (\vr);
    \draw (cwi) [fill=white] circle (\vr);
    
    \node at (av) [above] {$a_v$};
    \node at (awi) [left] {$a_{w_i}$};
    \node at (bwi) [below] {$b_{w_i}$};
    \node at (cwi) [below] {$c_{w_i}$};
    
    \node at (1.0, -0.3) {$T_{w_i}$};
    \node at (3.85, -1) {$\dots$};
\end{scope}

\begin{scope}[shift={(11.0,0)}]
    \path (1.0, 1.0) coordinate (awl);
    \path (-0.5, -1.0) coordinate (bwl);
    \path (2.5, -1.0) coordinate (cwl);
    \path (-1.0, -1.0) coordinate (bwl1);
    \path (2.5, 3.0) coordinate (cwl1);

    \draw (awl) -- (bwl);
    \draw (awl) -- (cwl);
    \draw (bwl) -- (cwl);
    \draw (bwl) -- (bwl1);
    \draw (cwl) -- (cwl1);
    \draw[-latex, shorten >=.65mm, shorten <=3.5mm, line width=0.25mm] (bwl1) -- (bwl);
    \draw[-latex, shorten >=.65mm, shorten <=3.5mm, line width=0.25mm] (bwl) -- (cwl);
    \draw[-latex, shorten >=.65mm, shorten <=3.5mm, line width=0.25mm] (cwl) -- (cwl1);
    
    \draw (awl) [fill=blue] circle (\vr);
    \draw[dashed, red] (awl) circle (28pt);
    \draw (bwl) [fill=blue!20] circle (\vr);
    \draw (cwl) [fill=blue!20] circle (\vr);
    \draw (cwl1) [fill=blue!20] circle (\vr);
    \draw[dashed, red] (cwl1) circle (18pt);
    
    \node at (awl) [left] {$a_{w_{\ell}}$};
    \node at (bwl) [below] {$b_{w_{\ell}}$};
    \node at (cwl) [below] {$c_{w_{\ell}}$};
    \node at (cwl1) [above] {$a_y$};
    
    \node at (1.0, -0.3) {$T_{w_{\ell}}$};
\end{scope}

\end{tikzpicture}
\end{center}
\caption{An illustration of applying Extension Rule~\ref{swap:paths}.}
\label{fig:path-fix-2}
\end{figure}

Let the zero forcing process begin from $S'$ until $a_{w_i}$, after which, $c_{w_i}$ would then be the only white-colored neighbor of $a_{w_i}$. Thus, $a_{w_i}$ may then force $c_{w_i}$ to become blue colored. After $c_{w_i}$ becomes blue-colored, $w_i$ would be a $\frac{3}{3}$-blue-colored vertex of $M_G$. Moreover, $w_i$ is connected to the $\frac{1}{3}$-blue-colored neighbors $w_{i-1}$ and $w_{i+1}$ by blue-white-open edges $c_{w_{i-1}}b_{w_i}$ and $c_{w_i}b_{w_{i+1}}$, respectively. By Claim~\ref{claim:1-blue-one}, $w_{i-1}$ and $w_{i+1}$ would then become $\frac{3}{3}$-blue-colored. Thereafter, this process would continue until both $w_1$ and $w_{\ell}$ become $\frac{3}{3}$-blue-colored. After $w_{\ell}$ becomes $\frac{3}{3}$-blue-colored, the only white colored neighbor of $c_{w_{\ell}}$ would be $a_{y''}$, and so, $c_{w_{\ell}}$ may then force $a_{y''}$ to become blue-colored; see Figure~\ref{fig:path-fix-2} for an illustration. Note that all vertices in the dominating set $X$ of $G$ remain reachable from $S'$. Thus, $N_G^z[S] \subset N_G^z[S']$, and the proof of the claim is finished. \smallqed

If $H_G^w$ was not a path component of $G[Y]$, then $H_G^w$ is necessarily a cycle component of $G[Y]$. Denote to the vertices of this component by $H_G^w: b_{w_1}c_{w_1}\dots b_{w_{\ell}}c_{w_{\ell}}b_{w_1}$, and denote the associated vertices (triangle-units) in $M_G$ by $H_M^w: w_1w_2\dots w_{\ell}w_1$. Note for some $i \in [\ell]$, $a_{w_i} = a_w$, $b_{w_i} = b_w$, and $c_{w_i} = c_w$. configuration.

\begin{figure}[htb]
 \begin{center}
\begin{tikzpicture}[scale=.8,style=thick,x=1cm,y=1cm]
\def\vr{3.5pt}

\path (1.0, 1.0) coordinate (aw1);
\path (1.0, 2.5) coordinate (aw11);
\path (-0.5, -1.0) coordinate (bw1);
\path (-0.5, -3.0) coordinate (bw11);
\path (2.5, -1.0) coordinate (cw1);
\path (3, -1.0) coordinate (cw11);

\draw (aw1) -- (bw1);
\draw (aw1) -- (aw11);
\draw (aw1) -- (cw1);
\draw (bw1) -- (cw1);
\draw[bend right=105] (bw1) to (bw11);
\draw (cw1) -- (cw11);

\draw (aw1) [fill=blue] circle (\vr);
\draw[dashed, red] (aw1) circle (28pt);
\draw (bw1) [fill=white] circle (\vr);
\draw (cw1) [fill=white] circle (\vr);

\node at (aw1) [left] {$a_{w_1}$};
\node at (bw1) [below] {$b_{w_1}$};
\node at (cw1) [below] {$c_{w_1}$};

\node at (1.0, -0.3) {$T_{w_1}$};

\node at (3.85, -1) {$\dots$};

\begin{scope}[shift={(5.5,0)}]
    \path (1.0, 3.0) coordinate (av);
    \path (1.0, 1.0) coordinate (awi);
    \path (-0.5, -1.0) coordinate (bwi);
    \path (2.5, -1.0) coordinate (cwi);
    \path (3, -1.0) coordinate (cwi1);

    \path (-1.0, -1.0) coordinate (bwi1);

    \draw (awi) -- (bwi);
     \draw (bwi) -- (bwi1);
    \draw (awi) -- (cwi);
    \draw (bwi) -- (cwi);
    \draw (av) -- (awi);
    \draw (cwi) -- (cwi1);
    \draw[-latex, shorten >=.65mm, shorten <=3.5mm, line width=0.25mm] (av) -- (awi);
    
    \draw (av) [fill=blue] circle (\vr);
    \draw (awi) [fill=blue!20] circle (\vr);
    \draw[dashed, red] (awi) circle (28pt);
    \draw (bwi) [fill=white] circle (\vr);
    \draw (cwi) [fill=white] circle (\vr);
    
    \node at (av) [above] {$a_v$};
    \node at (awi) [left] {$a_{w_i}$};
    \node at (bwi) [below] {$b_{w_i}$};
    \node at (cwi) [below] {$c_{w_i}$};
    
    \node at (1.0, -0.3) {$T_{w_i}$};
    \node at (3.85, -1) {$\dots$};
\end{scope}

\begin{scope}[shift={(11.0,0)}]
    \path (1.0, 2.5) coordinate (awl1);
    \path (1.0, 1.0) coordinate (awl);
    \path (-0.5, -1.0) coordinate (bwl);
    \path (2.5, -1.0) coordinate (cwl);
    \path (-1.0, -1.0) coordinate (bwl1);
    \path (2.5, -3.0) coordinate (cwl1);


    
    \draw (awl) -- (awl1);
    \draw (awl) -- (bwl);
    \draw (awl) -- (cwl);
    \draw (bwl) -- (cwl);
    \draw (bwl) -- (bwl1);
    \draw[bend left=65] (cwl) to (cwl1);
    \draw (bw11) -- (cwl1);
    
    \draw (awl) [fill=blue] circle (\vr);
    \draw[dashed, red] (awl) circle (28pt);
    \draw (bwl) [fill=white] circle (\vr);
    \draw (cwl) [fill=white] circle (\vr);
    
    \node at (awl) [left] {$a_{w_{\ell}}$};
    \node at (bwl) [below] {$b_{w_{\ell}}$};
    \node at (cwl) [below] {$c_{w_{\ell}}$};
    
    \node at (1.0, -0.3) {$T_{w_{\ell}}$};
\end{scope}

\end{tikzpicture}
\end{center}
\caption{A dominating vertex $a_{w_i}$ is forced to become blue-colored and is on a cycle of the multigraph $M_G$.}
\label{fig:cycle-fix-1}
\end{figure}
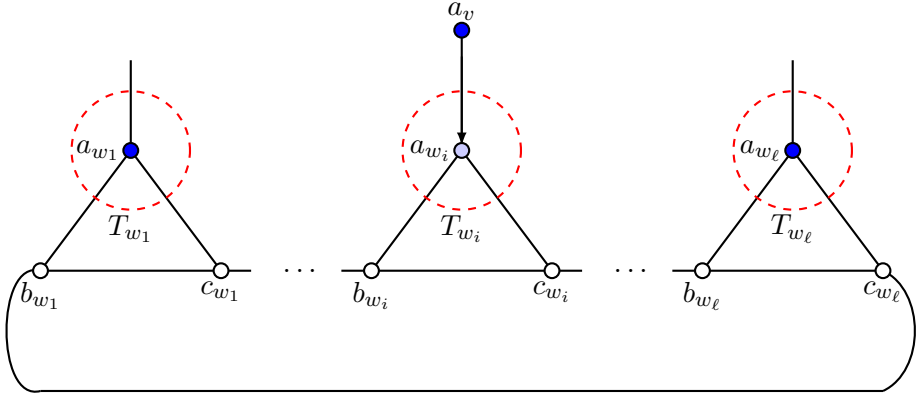

We introduce the following extension rule to deal with cycle components of $G[Y]$. 
\begin{swap}\label{swap:cycles}
Let $a_v$ be a vertex that forced the dominating vertex $a_{w_i}$ in the triangle-unit $T_{w_i}$ with vertex set $V(T_{w_i}) = \{a_{w_i}, b_{w_i}, c_{w_i} \}$. If the component $H_G$ containing $b_{w_i}$ and $c_{w_i}$ is a cycle, then
\[
S \leftarrow (S \setminus \{a_{w_{\ell}}\}) \cup \{b_{w_i}\}, 
\]
where $a_{y''}$ is adjacent with $c_{w_{\ell}}$ (and $a_{y''} \notin T_{w_i}$). 
\end{swap}

The following claim verifies the validity of Extension Rule~\ref{swap:cycles}.
\begin{claim}\label{claim:cycle-extention}
Extension Rule~\ref{swap:cycles} is a valid-$z$-proper-extension rule on $S$. 
\end{claim}
\proof Let $S' = (S \setminus \{a_{w_{\ell}}\}) \cup \{b_{w_i}\}$, and so, $|S'| = |S|$. Note that we have only swapped one vertex from the dominating set $X$ of $G$; see Extension Rule~\ref{swap:empty-triangle}. Thus, no vertex $a_{w_j}$, for $j \in [\ell] \setminus \{i\}$, has been forced to become blue-colored by the zero forcing process in $G$; that is, each of these vertices is initially blue-colored in $S$. If vertex adjacent with $a_{w_{\ell}}$ was contained in a $\frac{3}{3}$-blue-colored triangle-unit, then we could have applied Extension Rule~\ref{swap:easy} and then forced $w_{\ell}$ to become $\frac{3}{3}$-blue-colored. If $a_{w_{\ell}}$ was not contained in a $\frac{3}{3}$-blue-colored triangle-unit, then necessarily no neighbor of $a_{w_{\ell}}$ forced a color change during the zero forcing process starting with $S$. Thus, all color changes under the set $S$ will also occur from starting with set $S'$. 

\begin{figure}[htb]
 \begin{center}
\begin{tikzpicture}[scale=.8,style=thick,x=1cm,y=1cm]
\def\vr{3.5pt}

\path (1.0, 1.0) coordinate (aw1);
\path (1.0, 2.5) coordinate (aw11);
\path (-0.5, -1.0) coordinate (bw1);
\path (-0.5, -3.0) coordinate (bw11);
\path (2.5, -1.0) coordinate (cw1);
\path (3, -1.0) coordinate (cw11);

\draw (aw1) -- (bw1);
\draw (aw1) -- (aw11);
\draw (aw1) -- (cw1);
\draw (bw1) -- (cw1);
\draw[bend right=105] (bw1) to (bw11);
\draw (cw1) -- (cw11);
\draw[-latex, shorten >=.65mm, shorten <=3.5mm, line width=0.25mm] (cw11) -- (cw1);
\draw[-latex, shorten >=.65mm, shorten <=3.5mm, line width=0.25mm] (cw1) -- (bw1);

\draw (aw1) [fill=blue] circle (\vr);
\draw[dashed, red] (aw1) circle (28pt);
\draw (bw1) [fill=blue!20] circle (\vr);
\draw (cw1) [fill=blue!20] circle (\vr);

\node at (aw1) [left] {$a_{w_1}$};
\node at (bw1) [below] {$b_{w_1}$};
\node at (cw1) [below] {$c_{w_1}$};

\node at (1.0, -0.3) {$T_{w_1}$};

\node at (3.85, -1) {$\dots$};

\begin{scope}[shift={(5.5,0)}]
    \path (1.0, 3.0) coordinate (av);
    \path (1.0, 1.0) coordinate (awi);
    \path (-0.5, -1.0) coordinate (bwi);
    \path (2.5, -1.0) coordinate (cwi);
    \path (3, -1.0) coordinate (cwi1);

    \path (-1.0, -1.0) coordinate (bwi1);

    \draw (awi) -- (bwi);
     \draw (bwi) -- (bwi1);
    \draw (awi) -- (cwi);
    \draw (bwi) -- (cwi);
    \draw (av) -- (awi);
    \draw (cwi) -- (cwi1);
    \draw[-latex, shorten >=.65mm, shorten <=3.5mm, line width=0.25mm] (bwi) -- (bwi1);
    \draw[-latex, shorten >=.65mm, shorten <=3.5mm, line width=0.25mm] (cwi) -- (cwi1);
    \draw[-latex, shorten >=.65mm, shorten <=3.5mm, line width=0.25mm] (av) -- (awi);
    \draw[-latex, shorten >=.65mm, shorten <=3.5mm, line width=0.25mm] (awi) -- (cwi);
    
    \draw (av) [fill=blue] circle (\vr);
    \draw (awi) [fill=blue!20] circle (\vr);
    \draw[dashed, red] (awi) circle (28pt);
    \draw (bwi) [fill=blue] circle (\vr);
    \draw (cwi) [fill=blue!20] circle (\vr);
    
    \node at (av) [above] {$a_v$};
    \node at (awi) [left] {$a_{w_i}$};
    \node at (bwi) [below] {$b_{w_i}$};
    \node at (cwi) [below] {$c_{w_i}$};
    
    \node at (1.0, -0.3) {$T_{w_i}$};
    \node at (3.85, -1) {$\dots$};
\end{scope}

\begin{scope}[shift={(11.0,0)}]
    \path (1.0, 2.5) coordinate (awl1);
    \path (1.0, 1.0) coordinate (awl);
    \path (-0.5, -1.0) coordinate (bwl);
    \path (2.5, -1.0) coordinate (cwl);
    \path (-1.0, -1.0) coordinate (bwl1);
    \path (2.5, -3.0) coordinate (cwl1);

    \path (1.0, 2.5) coordinate (aq);
    \path (-0.5, 4.0) coordinate (bq);
    \path (2.5, 4.0) coordinate (cq);

    
    \draw (awl) -- (awl1);
    \draw (awl) -- (bwl);
    \draw (awl) -- (cwl);
    \draw (bwl) -- (cwl);
    \draw (bwl) -- (bwl1);
    \draw[bend left=65] (cwl) to (cwl1);
    \draw (bw11) -- (cwl1);
    \draw[-latex, shorten >=.65mm, shorten <=3.5mm, line width=0.25mm] (bwl1) -- (bwl);
    \draw[-latex, shorten >=.65mm, shorten <=3.5mm, line width=0.25mm] (bw11) -- (cwl1);
    \draw[-latex, shorten >=.65mm, shorten <=3.5mm, line width=0.25mm] (cwl) -- (awl);
    
    \draw (awl) [fill=blue!20] circle (\vr);
    \draw[dashed, red] (awl) circle (28pt);
    \draw (bwl) [fill=blue!20] circle (\vr);
    \draw (cwl) [fill=blue!20] circle (\vr);

    
    \node at (awl) [left] {$a_{w_{\ell}}$};
    \node at (bwl) [below] {$b_{w_{\ell}}$};
    \node at (cwl) [below] {$c_{w_{\ell}}$};

    
    \node at (1.0, -0.3) {$T_{w_{\ell}}$};
\end{scope}

\end{tikzpicture}
\end{center}
\caption{An illustration of applying Extension Rule~\ref{swap:cycles}.}
\label{fig:cycle-fix-2}
\end{figure}
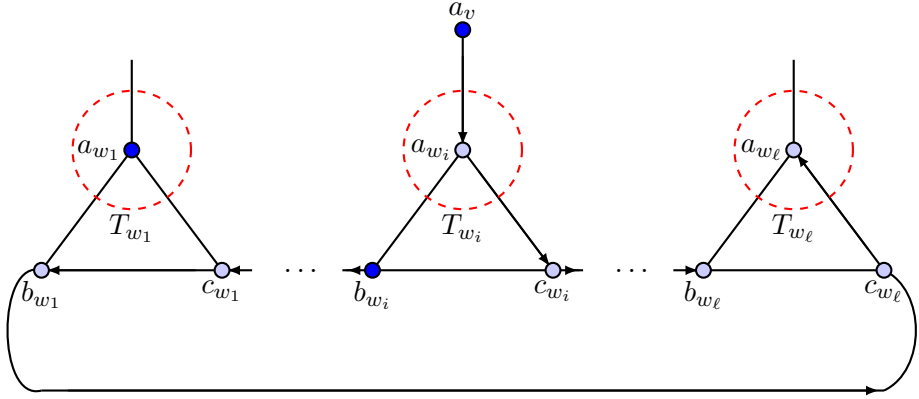

Starting from the set $S'$, let the zero forcing process start and continue. At some iteration of the zero forcing process, $a_v$ will force $a_{w_i}$ to become blue-colored. Thereafter, $c_{w_i}$ would be the only white-colored neighbor of $a_{w_i}$. Thus, $a_{w_i}$ may then force $c_{w_i}$ to become blue-colored. After $c_{w_i}$ becomes blue-colored, $w_i$ would be $\frac{3}{3}$-blue-colored. Thus, $w_i$ may then force onto the triangle-units $T_{w_{i-1}}$ and $T_{w_{i+1}}$ along the blue-white open edges $b_{w_i}c_{w_{i-1}}$ and $c_{w_i}b_{w_{i+1}}$, respectively. In particular, a sequence of color changes occurs until $w_1$ becomes $\frac{3}{3}$-blue-colored. Thereafter, $c_{w_{\ell}}$ is the only white-colored neighbor of $b_{w_1}$, and so, $b_{w_1}$ may then force $c_{w_{\ell}}$ to become blue-colored. Note that $w_{\ell}$ was also forced onto by $w_{\ell-1}$, and so, $b_{w_{\ell}}$ is blue-colored. Thus, at some iteration of the zero forcing process starting with $S'$, $a_{w_{\ell}}$ would be the only white-colored neighbor of $c_{w_{\ell}}$. Hence, $c_{w_{\ell}}$ may then force $a_{w_{\ell}}$ to become colored blue; see Figure~\ref{fig:cycle-fix-2} for an illustration. Therefore, $N_G^z[S] \subset N_G^z[S']$, and the proof of the claim is finished. \smallqed

Apply Extension Rule~\ref{swap:paths} if $H_G^w$ is a path component of $G[Y]$, and Extension Rule~\ref{swap:cycles} if $H_G^w$ is a cycle component of $G[Y]$. Thereafter, allow the zero forcing process to continue and apply Extension Rule~\ref{swap:easy} whenever possibly until no further color changes are possible. If a triangle-unit that we swapped a dominating vertex does not become $\frac{3}{3}$-blue-colored, then apply either Extension Rule~\ref{swap:paths} or Extension Rule~\ref{swap:cycles}, depending on the component of $G[Y]$ in consideration. After all, borrowed from triangle-units become $\frac{3}{3}$-blue-colored, apply Extension Rule~\ref{swap:empty-triangle} if need be and repeat. Since $G$ is connected, we eventually arrive at a set $S$ of blue-colored vertices so that all of $V(G)$ becomes blue-colored, and so, the resulting set $S$ is a zero forcing set of $G$. Thus, 
\[
Z(G) \leq |S| = |X| + 1 = \gamma(G) + 1, 
\]
and the proof of the lemma is finished. 
\qed

In the next section, we give our proof for Conjecture~\ref{conj:confirmed}. 

\section{Proof of Conjecture~\ref{conj:confirmed}}
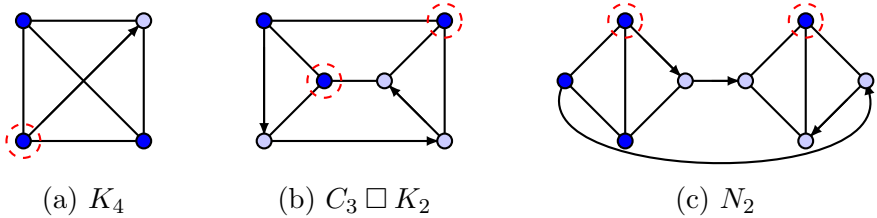
\begin{figure}[htb]
\begin{center}
\begin{tikzpicture}[scale=.8,style=thick,x=1cm,y=1cm]
\def\vr{3.5pt}

\path (-5.0, 2.0) coordinate (v1);
\path (-5.0, 0.0) coordinate (v2);
\path (-3, 2.0) coordinate (v3);
\path (-3, 0.0) coordinate (v4);

\draw (v1) -- (v2);
\draw (v1) -- (v3);
\draw (v1) -- (v4);
\draw (v2) -- (v3);
\draw (v2) -- (v4);
\draw (v3) -- (v4);

\draw[-latex, shorten >=.65mm, shorten <=3.5mm, line width=0.25mm] (v2) -- (v3);

\draw (v1) [fill=blue] circle (\vr);
\draw (v2) [fill=blue] circle (\vr);
\draw (v3) [fill=blue!20] circle (\vr);
\draw (v4) [fill=blue] circle (\vr);
\draw[dashed, red] (v2) circle (8pt);

\node at (-4.0, -1) {(a) $K_4$};
\path (0, 1) coordinate (e);
\path (-1.0, 2.0) coordinate (e1);
\path (-1.0, 0.0) coordinate (e2);

\path (1, 1) coordinate (f);
\path (2, 2.0) coordinate (f1);
\path (2, 0.0) coordinate (f2);

\draw (e) -- (e1);
\draw (e) -- (e2);
\draw (e1) -- (e2);

\draw (f) -- (f1);
\draw (f) -- (f2);
\draw (f1) -- (f2);

\draw (e) -- (f);
\draw (e1) -- (f1);
\draw (e2) -- (f2);

\draw[-latex, shorten >=.65mm, shorten <=3.5mm, line width=0.25mm] (e1) -- (e2);
\draw[-latex, shorten >=.65mm, shorten <=3.5mm, line width=0.25mm] (e2) -- (f2);
\draw[-latex, shorten >=.65mm, shorten <=3.5mm, line width=0.25mm] (f2) -- (f);

\draw (e) [fill=blue] circle (\vr);
\draw (e1) [fill=blue] circle (\vr);
\draw (e2) [fill=blue!20] circle (\vr);

\draw[dashed, red] (e) circle (8pt);
\draw[dashed, red] (f1) circle (8pt);

\draw (f) [fill=blue!20] circle (\vr);
\draw (f1) [fill=blue] circle (\vr);
\draw (f2) [fill=blue!20] circle (\vr);

\node at (0.5, -1) {(b) $C_3 \: \Box \: K_2$};

\path (4, 1) coordinate (d1);
\path (5, 2.0) coordinate (d2);
\path (5, 0.0) coordinate (d3);
\path (6, 1.0) coordinate (d4);

\path (7, 1) coordinate (dd1);
\path (8, 2.0) coordinate (dd2);
\path (8, 0.0) coordinate (dd3);
\path (9, 1.0) coordinate (dd4);

\draw (d1) -- (d2);
\draw (d1) -- (d3);
\draw (d2) -- (d3);
\draw (d2) -- (d4);
\draw (d3) -- (d4);

\draw (d4) -- (dd1);

\draw[-latex, shorten >=.65mm, shorten <=3.5mm, line width=0.25mm] (d2) -- (d4);
\draw[-latex, shorten >=.65mm, shorten <=3.5mm, line width=0.25mm] (d4) -- (dd1);
\draw[-latex, shorten >=.65mm, shorten <=3.5mm, line width=0.25mm] (dd4) -- (dd3);
\draw[->, bend right, >=latex] (d1) to[out=-120, in=-75] (dd4);

\draw (dd1) -- (dd2);
\draw (dd1) -- (dd3);
\draw (dd2) -- (dd3);
\draw (dd2) -- (dd4);
\draw (dd3) -- (dd4);

\draw (d1) [fill=blue] circle (\vr);
\draw (d2) [fill=blue] circle (\vr);
\draw (d3) [fill=blue] circle (\vr);
\draw (d4) [fill=blue!20] circle (\vr);
\draw (dd1) [fill=blue!20] circle (\vr);
\draw (dd2) [fill=blue] circle (\vr);
\draw (dd3) [fill=blue!20] circle (\vr);
\draw (dd4) [fill=blue!20] circle (\vr);

\draw[dashed, red] (d2) circle (8pt);
\draw[dashed, red] (dd2) circle (8pt);

\node at (6.5, -1) {(c) $N_2$};

\end{tikzpicture}
\end{center}
\caption{The complete graph $K_4$, the prism $C_3 \: \Box \: K_2$, and the diamond-necklace $N_2$ as seen in the base cases for the proof of Theorem~\ref{conj:confirmed}. Minimum zero forcing sets are shown in dark blue, while forcing steps and color changes are indicated by directed edges and light blue colored vertices, respectively; minimum dominating sets shown by red dashed circles.}
\label{fig:main-thm-0}
\end{figure}

In this section, we prove Conjecture~\ref{conj:confirmed} and characterize graphs attaining equality in its statement. 
\begin{thm}\label{thm:main}
If $G$ is a connected, cubic, and claw-free graph, then 
\[
Z(G) \leq \gamma(G) + 2, 
\]
with equality if and only if $G \in \mathcal{N}_{\text{cubic}}^*$. 
\end{thm}
\proof We will proceed by induction on the order $n\geq 4$ of a connected, cubic, and claw-free graph. If $n = 4$, then $G = K_4$ implying $Z(G) = \gamma(G) + 2$; see Figure~\ref{fig:main-thm-0} (a). If $n = 6$, then $G$ is the prism $C_3 \: \Box \: K_2$ implying $Z(G) = 3 < \gamma(G) + 2 = 4$; see Figure~\ref{fig:main-thm-0} (b). If $n = 8$, then $G$ is the diamond-necklace $N_2$ which implies $Z(G) = \gamma(G) + 2$ by Lemma~\ref{lem:diamond-formula}; also see Figure~\ref{fig:main-thm-0} (c). Thus, we have established our base cases. Next, let $n \geq 10$ and assume that if $G'$ is a connected, cubic, and claw-free graph of order $n'$, where $6 \leq n' < n$, then $Z(G') \leq \gamma(G') + 2$.

Let $G$ be a connected, cubic, and claw-free graph with order $n$. If $G \in \mathcal{N}_{\text{cubic}}^*$, then by Lemma~\ref{lem:diamond-formula}, $Z(G) = \gamma(G) + 2$. Thus, we may assume that $G \notin \mathcal{N}_{\text{cubic}}^*$, since otherwise the desired result follows. Hence, at least one unit in the $\Delta$-D-partition of $V(G)$ is a diamond-unit. Since every triangle-unit of $G$ is joined by three edges to vertices in other units, and since every diamond-unit is joined by two edges to vertices in other units, it must be the case that there are at least two triangle-units in our $\Delta$-D-partition. If $V(G)$ does not contain a diamond-unit in its $\Delta$-D-partition, then $G$ has a spanning 2-factor consisting of triangles, which by Lemma~\ref{lem:2-factor}, implies $Z(G) \leq \gamma(G) + 1 < \gamma(G) + 2$. Thus, we may further assume $G$ contains at least one diamond since otherwise, the desired result follows. Let $D$ be an arbitrary diamond-unit in the $\Delta$-D-partition of $V(G)$ and denote the vertices of $D$ by $V(D) = \{a, b, c, d\}$, where $ab$ is the missing edge in $D$. Let $e$ be the neighbor of $a$, not in $D$, and $f$ be the neighbor of $b$, not in $f$, where $e \neq f$ because $G$ is claw-free. The vertices $e$ and $f$ may, or may not be adjacent, and we consider these cases with the following claims. 

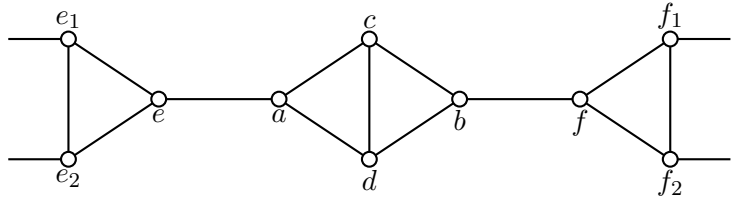
\begin{figure}[htb]
\begin{center}
\begin{tikzpicture}[scale=.8,style=thick,x=1cm,y=1cm]
\def\vr{3.5pt}


\path (-4.0, 1.0) coordinate (e1);
\path (-5.0, 1.0) coordinate (e11);
\path (-4.0, -1.0) coordinate (e2);
\path (-5.0, -1.0) coordinate (e22);
\path (-2.5, 0.0) coordinate (e);

\path (-0.5, 0.0) coordinate (a);
\path (1.0, 1.0) coordinate (c);
\path (1.0, -1.0) coordinate (d);
\path (2.5, 0.0) coordinate (b);

\path (4.5, 0.0) coordinate (v1);
\path (6, 1.0) coordinate (v2);
\path (7, 1.0) coordinate (v22);
\path (6, -1.0) coordinate (v3);
\path (7, -1.0) coordinate (v33);

\draw (e) -- (e1);
\draw (e) -- (e2);
\draw (e1) -- (e2);
\draw (e1) -- (e11);
\draw (e2) -- (e22);
\draw (e) -- (a);

\draw (a) -- (c);
\draw (a) -- (d);
\draw (c) -- (d);
\draw (c) -- (b);
\draw (d) -- (b);

\draw (b) -- (v1);
\draw (v1) -- (v2);
\draw (v1) -- (v3);
\draw (v2) -- (v3);

\draw (v2) -- (v22);
\draw (v3) -- (v33);

\draw (e) [fill=white] circle (\vr);
\draw (e1) [fill=white] circle (\vr);
\draw (e2) [fill=white] circle (\vr);

\draw (a) [fill=white] circle (\vr);
\draw (b) [fill=white] circle (\vr);
\draw (c) [fill=white] circle (\vr);
\draw (d) [fill=white] circle (\vr);

\draw (v1) [fill=white] circle (\vr);
\draw (v2) [fill=white] circle (\vr);
\draw (v3) [fill=white] circle (\vr);

\node at (e) [below] {$e$};
\node at (e1) [above] {$e_1$};
\node at (e2) [below] {$e_2$};

\node at (a) [below] {$a$};
\node at (b) [below] {$b$};
\node at (c) [above] {$c$};
\node at (d) [below] {$d$};

\node at (v1) [below] {$f$};
\node at (v2) [above] {$f_1$};
\node at (v3) [below] {$f_2$};

\end{tikzpicture}
\end{center}
\caption{The structure of the subgraph as seen in the proof of Claim~\ref{claim:thm-1}. }
\label{fig:main-thm-1}
\end{figure}

\begin{claim}\label{claim:thm-1}
If $e$ and $f$ are not adjacent, then $Z(G) <  \gamma(G) + 2$. 
\end{claim}
\proof Suppose that $e$ and $f$ are not adjacent. Next, let $e_1$ and $e_2$ be the neighbors of $e$ different from $a$, and let $f_1$ and $f_2$ be the neighbors of $f$ different from $b$. Since $G$ is claw-free, note that $\{e, e_1, e_2\}$ and $\{f, f_1, f_2\}$ both induce triangles in $G$, say $T_e$ and $T_f$, respectively. Note that $T_e$ and $T_f$ share no common vertex. Moreover, $T_e$ and $T_f$ are themselves either triangle-units, or are a part of some diamond-unit in the $\Delta$-D-partition of $V(G)$; see Figure~\ref{fig:main-thm-1} for an illustration of this configuration.

\begin{figure}[htb]
\begin{center}
\begin{tikzpicture}[scale=.8,style=thick,x=1cm,y=1cm]
\def\vr{3.5pt}


\path (-4.0, 2.0) coordinate (e1);
\path (-4.0, -1.5) coordinate (e2);
\path (-5.0, -1.0) coordinate (e22);
\path (-2.5, 0.25) coordinate (e);

\path (-0.5, 0.25) coordinate (a);
\path (1.0, 1.25) coordinate (c);
\path (1.0, -0.75) coordinate (d);
\path (2.5, 0.25) coordinate (b);

\path (4.5, 0.25) coordinate (f);
\path (6, 2.0) coordinate (f1);
\path (6, -1.5) coordinate (f2);

\draw (e) -- (e1);
\draw (e) -- (e2);
\draw (e1) -- (e2);
\draw (e1) -- (f1);
\draw (e2) -- (f2);
\draw (e) -- (a);

\draw (a) -- (c);
\draw (a) -- (d);
\draw (c) -- (d);
\draw (c) -- (b);
\draw (d) -- (b);

\draw (b) -- (f);
\draw (f) -- (f1);
\draw (f) -- (f2);
\draw (f1) -- (f2);

\draw[-latex, shorten >=.65mm, shorten <=3.5mm, line width=0.25mm] (e1) -- (e2);
\draw[-latex, shorten >=.65mm, shorten <=3.5mm, line width=0.25mm] (e2) -- (f2);
\draw[-latex, shorten >=.65mm, shorten <=3.5mm, line width=0.25mm] (f2) -- (f);
\draw[-latex, shorten >=.65mm, shorten <=3.5mm, line width=0.25mm] (f) -- (b);
\draw[-latex, shorten >=.65mm, shorten <=3.5mm, line width=0.25mm] (b) -- (d);
\draw[-latex, shorten >=.65mm, shorten <=3.5mm, line width=0.25mm] (d) -- (a);

\draw (e) [fill=blue] circle (\vr);
\draw[dashed, red] (e) circle (18pt);
\draw (e1) [fill=blue] circle (\vr);
\draw (e2) [fill=blue!20] circle (\vr);

\draw (a) [fill=blue!20] circle (\vr);
\draw (b) [fill=blue!20] circle (\vr);
\draw (c) [fill=blue] circle (\vr);
\draw[dashed, red] (c) circle (18pt);
\draw (d) [fill=blue!20] circle (\vr);

\draw (f) [fill=blue!20] circle (\vr);
\draw (f1) [fill=blue] circle (\vr);
\draw (f2) [fill=blue!20] circle (\vr);
\draw[dashed, red] (f) circle (18pt);

\node at (e) [below] {$e$};
\node at (e1) [above] {$e_1$};
\node at (e2) [below] {$e_2$};

\node at (a) [below] {$a$};
\node at (b) [below] {$b$};
\node at (c) [above] {$c$};
\node at (d) [below] {$d$};

\node at (f) [below] {$f$};
\node at (f1) [above] {$f_1$};
\node at (f2) [below] {$f_2$};

\end{tikzpicture}
\end{center}
\caption{The graph $G$ as seen in Claim~\ref{claim:thm-1} given in the proof of Theorem~\ref{conj:confirmed}. A minimum zero forcing set shown in dark blue, while forcing steps and color changes are indicated by directed edges and light blue colored vertices, respectively; a minimum dominating set shown by red dashed circles.}
\label{fig:main-thm-2}
\end{figure}
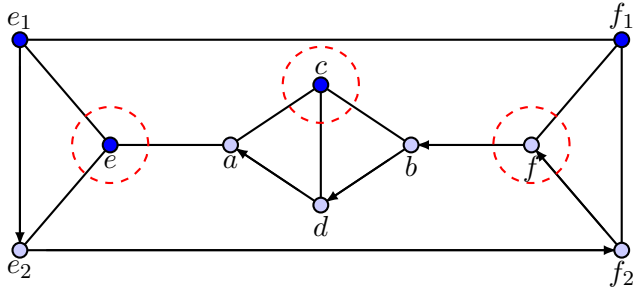

Let $G'$ be the graph obtained from $G$ by deleting the vertices in $V(G)$ and any edges incident edges from $V(G)$, and then adding the edge $ef$. Note that $G'$ is a connected, cubic, and claw-free graph with order $n'$, where $6 \leq n' < n$. If the $\Delta$-D-partition of $V(G')$ consists of only diamond-units and no triangle-units, then $G \in \mathcal{N}_{\text{cubic}}^*$, a contradiction our earlier assumption that $G \notin \mathcal{N}_{\text{cubic}}^*$. Thus, $G' \notin \mathcal{N}_{\text{cubic}}$. Hence, the $\Delta$-D-partition of $V(G')$ contains at least one triangle-unit. If $G'$ is the prism $C_3 \: \Box \: K_2$, then renaming vertices if necessary, we may assume $e_i$ and $f_i$ are adjacent for $i \in [2]$, which implies that $G$ is the graph with order $n = 10$ shown in Figure~\ref{fig:main-thm-2}. The set $\{c, e, f_1\}$ is a minimum dominating set of $G$, and the set $\{c, e, e_1, f \}$ is a minimum zero forcing set of $G$. Thus, $Z(G) = 4 < \gamma(G) + 2 = 5$. Hence, we may assume that $G'$ is not the prism $C_3 \: \Box \: K_2$ since otherwise, our desired result follows. We proceed by proving the following two subclaims.

\begin{subclaim}\label{subclaim:zero-forcing}
$Z(G) = Z(G') + 1$
\end{subclaim}
\proof Let $G$ be a cubic and claw-free graph. Next, let $D$ be an arbitrary diamond-unit in the $\Delta$-D-partition of $V(G)$ and denote the vertices of $D$ by $V(D) = \{a, b, c, d\}$, where $ab$ is the missing edge in $D$. Let $e$ be the neighbor of $a$, not in $D$, and $f$ be the neighbor of $b$, not in $f$, where $e \neq f$ because $G$ is claw-free.Let $S' \subseteq V(G')$ be a minimum zero forcing set of $G'$. If $f \notin S'$, then at some iteration of the zero forcing process in $G'$, a neighbor of $f$ forces $f$ to become colored blue. If $f$ was forced to become colored blue during the zero forcing process in $G'$, let $S = S' \cup \{c\}$ in $G$. If $e$ was the vertex that forced $f$ to become colored blue in $G'$, starting from the blue coloring $S$ in $G$, at some iteration of the zero forcing process $e$ would force $a$ to become colored blue. After $a$ becomes colored blue, $d$ would be the only white-colored neighbor of $a$, so $a$ would then force $d$ to become colored blue. After $d$ becomes colored blue, $b$ would be the only white-colored neighbor of $d$, so $d$ would then force $b$ to become colored blue. Finally, after $b$ becomes colored blue, $f$ would be the only white-colored neighbor of $f$, so $d$ would force $f$ to become colored blue. After $f$ becomes colored blue, all color changes in $G'$ starting from $S'$ would now occur in $G$. Thus, $S$ is a zero forcing set of $G$. Hence, $Z(G) \leq |S'| + 1 = Z(G') + 1$. However, since every diamond-unit contains at least one vertex from every zero forcing set of $G$, it must be the case that $Z(G) = Z(G') + 1$.

Next, suppose that $f$ was forced to become colored blue by either $f_1$ or $f_2$. Without loss of generality, suppose $f_2$ forces $f$ to become colored blue. At some iteration of the zero forcing process in $G$ $f$ is the only white-colored neighbor of $f_2$, so $f_1$ was colored blue before $f$ was colored blue. These color changes would still occur in $G$ starting from the set $S$, after which the only white-colored neighbor of $f$ would be $b$, and so, $f$ may then force $b$ to become colored blue. After $b$ becomes colored blue, $d$ would be the only white-colored neighbor, so $b$ would force $d$ to become colored blue. After $d$ becomes colored blue, $a$ would be the only white-colored neighbor of $d$, so $d$ would then force $a$ to become colored blue. At this point of the zero forcing process, all colors that happened in $G'$ would still happen in $G$. Thus, $S$ is a zero forcing set of $G$. Hence, $Z(G) \leq |S'| + 1 = Z(G') + 1$. However, since every diamond-unit contains at least one vertex from every zero forcing set of $G$, it must be the case that $Z(G) = Z(G') + 1$. Hence, the result follows if $f$ was forced to become colored blue during the zero forcing process on $G'$. 

Since the above arguments apply to the vertex $e$, we may assume that neither $e$ nor $f$ become colored blue during the zero forcing process in $G'$. In particular, this implies $e, f \in S'$. Next suppose that at some iteration of the zero forcing process in $G'$, $f_2$ becomes colored blue. After $f_2$ becomes colored blue, the only white-colored neighbor of $f_2$ would be $f_1$, and so, $f_2$ may then force $f_1$ to become colored blue. After these color changes have in $G'$ and $G'$, the only white-colored neighbor of $f$ in $G$ would be $b$, so, $f$ could then force $b$ to become colored blue. After $b$ becomes colored blue, $d$ would be the only white-colored neighbor, so $b$ would force $d$ to become colored blue. After $d$ becomes colored blue, $a$ would be the only white-colored neighbor of $d$, so $d$ would then force $a$ to become colored blue. At this point of the zero forcing process, all colors that happened in $G'$ would still happen in $G$. Thus, $S$ is a zero forcing set of $G$. Hence, $Z(G) \leq |S'| + 1 = Z(G') + 1$. However, since every diamond-unit contains at least one vertex from every zero forcing set of $G$, it must be the case that $Z(G) = Z(G') + 1$. Hence, $Z(G) = Z(G') + 1$ when any neighbor of $f$ is forced to become colored blue by a vertex different than $f$ during the zero forcing process on $G'$.

Since the above arguments on $f$ also apply to the vertex $e$, we may also assume that no neighbor of $e$ is forced to become colored blue by a vertex different than $e$. Our suppositions combined imply that $e$ and $f$ each have exactly one neighbor, which is not contained in $S'$, where these white-colored neighbors are forced to become colored blue at the start of the zero forcing process in $G'$. Without loss of generality, suppose $e_2$ and $f_2$ are initially blue colored in $S'$, so $e$ forces $e_1$ to become colored blue and $f$ forces $f_1$ to become colored blue. Next let $S'' = (S'\setminus \{e\}) \cup \{f_1\}$ be a different set of initially blue colored vertices in $G'$, and then let $S = S'' \cup \{c \}$ and note $|S''| = |S'| = Z(G)$. Under the new coloring $S$ in $G$, observe that the only white-colored neighbor of $f$ in $G$ is $b$, and the only white-colored neighbor of $f$ in $G'$ is $e$. Thus, $f$ may force $b$ to become colored blue in $G$ and $e$ to become colored blue in $G'$. Thus, $S''$ is a zero forcing set, implying $|S''| = Z(G')$. Next, observe that after $b$ becomes colored blue in $G$, the only white-colored neighbor of $b$ would be $d$, and so, $b$ would force $d$ to become colored blue. After $d$ becomes colored blue, $a$ would be the only white-colored neighbor of $d$, so $d$ would then force $a$ to become colored blue. After $a$ becomes colored blue, $e$ would be the only white-colored neighbor of $a$, so $a$ may force $e$ to become colored blue. After $e$ becomes colored blue in $G$, all color changes in $G'$ will now occur in $G$. Thus, $S$ is a zero forcing set of $G$. Hence, $Z(G) \leq |S'| + 1 = Z(G') + 1$. However, since every diamond-unit contains at least one vertex from every zero forcing set of $G$, it must be the case that $Z(G) = Z(G') + 1$. 

\begin{subclaim}\label{subclaim:domination}
$\gamma(G) = \gamma(G') + 1$
\end{subclaim}
\proof Let $X' \subseteq V(G')$ be a minimum (independent) dominating set of $G'$. If both $e$ and $f$ are not in $X'$, then clearly $X = X' \cup \{c\}$ is a minimum (independent) dominating set of $G$ implying $\gamma(G) = \gamma(G') + 1$. Thus, we will $X'$ contains either $e$ or $f$. Note that since $e$ and $f$ are adjacent in $G'$, and since $X'$ is independent, either $e \in X'$ and $f \notin X'$, or $e \notin X'$ and $f \in X'$. Without loss of generality, suppose $e \notin X'$ and $f \in X'$. Then, the set $X = X' \cup \{a\}$ is clearly a minimum dominating set of $G$. Thus, $\gamma(G) = \gamma(G') + 1$, completing out proof. 
\smallqed 

Apply the inductive hypothesis to the graph $G'$, and so, 
\[
Z(G') < \gamma(G') + 2.
\]
By Claim~\ref{subclaim:zero-forcing}, $Z(G) = Z(G') + 1$, and by Claim~\ref{subclaim:domination}, $\gamma(G) = \gamma(G') + 1$. Thus, 
\[
Z(G) - 1= Z(G') < \gamma(G') + 2 = (\gamma(G) - 1) + 2.
\]
Hence,
\[
Z(G) < \gamma(G) + 2,
\]
and the proof of our claim is finished. \smallqed 

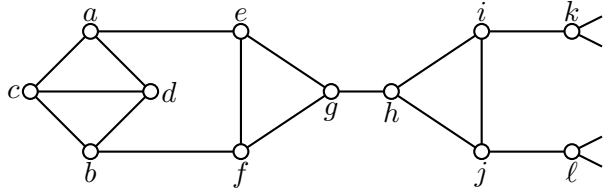
\begin{figure}[htb]
\begin{center}
\begin{tikzpicture}[scale=.8,style=thick,x=1cm,y=1cm]
\def\vr{3.5pt}

\path (1.0, 1.0) coordinate (a);
\path (0.0, 0.0) coordinate (c);
\path (1.0, -1.0) coordinate (b);
\path (2.0, 0.0) coordinate (d);

\path (3.5, 1.0) coordinate (e);
\path (3.5, -1.0) coordinate (f);
\path (5.0, 0.0) coordinate (g);

\path (7.5, 1.0) coordinate (i);
\path (7.5, -1.0) coordinate (j);
\path (6.0, 0.0) coordinate (h);

\path (9.0, 1.0) coordinate (k);
\path (9.0, -1.0) coordinate (ell);

\path (9.5, 1.25) coordinate (k1);
\path (9.5, 0.75) coordinate (k2);

\path (9.5, -0.75) coordinate (ell1);
\path (9.5, -1.25) coordinate (ell2);

\draw (a) -- (c);
\draw (a) -- (d);
\draw (c) -- (d);
\draw (c) -- (b);
\draw (d) -- (b);

\draw (a) -- (e);
\draw (b) -- (f);
\draw (e) -- (f);
\draw (e) -- (g);
\draw (f) -- (g);

\draw (g) -- (h);
\draw (h) -- (i);
\draw (h) -- (j);
\draw (i) -- (j);

\draw (i) -- (k);
\draw (j) -- (ell);

\draw (k) -- (k1);
\draw (k) -- (k2);

\draw (ell) -- (ell1);
\draw (ell) -- (ell2);

\draw (a) [fill=white] circle (\vr);
\draw (b) [fill=white] circle (\vr);
\draw (c) [fill=white] circle (\vr);
\draw (d) [fill=white] circle (\vr);
\draw (e) [fill=white] circle (\vr);
\draw (f) [fill=white] circle (\vr);
\draw (g) [fill=white] circle (\vr);
\draw (h) [fill=white] circle (\vr);
\draw (i) [fill=white] circle (\vr);
\draw (j) [fill=white] circle (\vr);
\draw (k) [fill=white] circle (\vr);
\draw (ell) [fill=white] circle (\vr);

\node at (a) [above] {$a$};
\node at (b) [below] {$b$};
\node at (c) [left] {$c$};
\node at (d) [right] {$d$};
\node at (e) [above] {$e$};
\node at (f) [below] {$f$};
\node at (g) [below] {$g$};
\node at (h) [below] {$h$};
\node at (i) [above] {$i$};
\node at (j) [below] {$j$};
\node at (k) [above] {$k$};
\node at (ell) [below] {$\ell$};

\end{tikzpicture}
\end{center}
\caption{The subgraph of $G$ appearing towards the end of the proof given for Theorem~\ref{conj:confirmed}.}
\label{fig:main-thm-3}
\end{figure}

By Claim~\ref{claim:thm-1}, we may assume that $e$ and $f$ are adjacent in $G$ since otherwise, the desired inequality follows. Thus, $e$ and $f$ belong to a common triangle-unit in the $\Delta$-D-partition of $V(G)$, say $T$. Let $g$ be the remaining vertex in the triangle-unit $T$, and let $h$ be the neighbor of $g$ not in $T$. Further, let $i$ and $j$ be the two vertices in the triangle containing $h$. If $h$ belongs to a diamond-unit, then choosing this diamond-unit initially in our argument would bring us back to Claim~\ref{claim:thm-1}, implying $Z(G) \leq \gamma(G) + 1$. Thus, we may assume that $h$ belongs to a triangle-unit in the $\Delta$-D-partition of $V(G)$; that is $\{h, i, j\}$ form a triangle-unit. Let $k$ and $\ell$ denote the neighbors of $i$ and $j$, respectively, not contained in $\{h, i, j\}$. By assumption, the vertex $h$ does not belong to a diamond-unit, so $k \neq \ell$. The resulting subgraph of $G$ is illustrated in Figure~\ref{fig:main-thm-3}, where $k$ and $\ell$ may be adjacent vertices.  

Suppose now that the vertex $k$ belongs to a diamond-unit $D^*$ in $G$, which implies that $D^*$ also is a diamond-unit of $G'$ considered previously. If the diamond-unit $D^*$ does not contain the vertex $\ell$, then choosing this diamond-unit as our initial diamond-unit $D$ brings us back to Claim~\ref{claim:thm-1}, which implies that $Z(G) \leq \gamma(G) + 2$. 

\begin{figure}[htb]
\begin{center}
\begin{tikzpicture}[scale=.8,style=thick,x=1cm,y=1cm]
\def\vr{3.5pt}

\path (1.0, 1.0) coordinate (a);
\path (0.0, 0.0) coordinate (c);
\path (1.0, -1.0) coordinate (b);
\path (2.0, 0.0) coordinate (d);

\path (3.5, 1.0) coordinate (e);
\path (3.5, -1.0) coordinate (f);
\path (5.0, 0.0) coordinate (g);

\path (7.5, 1.0) coordinate (i);
\path (7.5, -1.0) coordinate (j);
\path (6.0, 0.0) coordinate (h);

\path (10.0, 1.0) coordinate (k);
\path (10.0, -1.0) coordinate (ell);

\path (9.0, 0.0) coordinate (p);
\path (11.0, 0.0) coordinate (m);

\draw (a) -- (c);
\draw (a) -- (d);
\draw (c) -- (d);
\draw (c) -- (b);
\draw (d) -- (b);

\draw (a) -- (e);
\draw (b) -- (f);
\draw (e) -- (f);
\draw (e) -- (g);
\draw (f) -- (g);

\draw (g) -- (h);
\draw (h) -- (i);
\draw (h) -- (j);
\draw (i) -- (j);

\draw (i) -- (k);
\draw (j) -- (ell);

\draw (k) -- (p);
\draw (k) -- (m);
\draw (ell) -- (p);
\draw (ell) -- (m);
\draw (p) -- (m);

\draw[-latex, shorten >=.65mm, shorten <=3.5mm, line width=0.25mm] (a) -- (d);
\draw[-latex, shorten >=.65mm, shorten <=3.5mm, line width=0.25mm] (d) -- (b);
\draw[-latex, shorten >=.65mm, shorten <=3.5mm, line width=0.25mm] (b) -- (f);

\draw[-latex, shorten >=.65mm, shorten <=3.5mm, line width=0.25mm] (f) -- (g);

\draw[-latex, shorten >=.65mm, shorten <=3.5mm, line width=0.25mm] (g) -- (h);

\draw[-latex, shorten >=.65mm, shorten <=3.5mm, line width=0.25mm] (h) -- (j);

\draw[-latex, shorten >=.65mm, shorten <=3.5mm, line width=0.25mm] (j) -- (ell);

\draw[-latex, shorten >=.65mm, shorten <=3.5mm, line width=0.25mm] (ell) -- (p);

\draw[-latex, shorten >=.65mm, shorten <=3.5mm, line width=0.25mm] (p) -- (k);

\draw (a) [fill=blue] circle (\vr);
\draw (b) [fill=blue!20] circle (\vr);
\draw (c) [fill=blue] circle (\vr);
\draw (d) [fill=blue!20] circle (\vr);
\draw (e) [fill=blue] circle (\vr);
\draw (f) [fill=blue!20] circle (\vr);
\draw (g) [fill=blue!20] circle (\vr);
\draw (h) [fill=blue!20] circle (\vr);
\draw (i) [fill=blue] circle (\vr);
\draw (j) [fill=blue!20] circle (\vr);
\draw (k) [fill=blue!20] circle (\vr);
\draw (ell) [fill=blue!20] circle (\vr);
\draw (p) [fill=blue!20] circle (\vr);
\draw (m) [fill=blue] circle (\vr);

\node at (a) [above] {$a$};
\node at (b) [below] {$b$};
\node at (c) [left] {$c$};
\node at (d) [right] {$d$};
\node at (e) [above] {$e$};
\node at (f) [below] {$f$};
\node at (g) [below] {$g$};
\node at (h) [below] {$h$};
\node at (i) [above] {$i$};
\node at (j) [below] {$j$};
\node at (k) [above] {$k$};
\node at (ell) [below] {$\ell$};
\node at (p) [left] {$p$};
\node at (m) [right] {$m$};

\end{tikzpicture}
\end{center}
\caption{The graph $G$ appearing in the proof of Theorem~\ref{conj:confirmed}. A minimum zero forcing set is shown in dark blue, whereas zero forcing color changes are indicated by directed edges and lighter blue colored vertices. }
\label{fig:main-thm-4}
\end{figure}
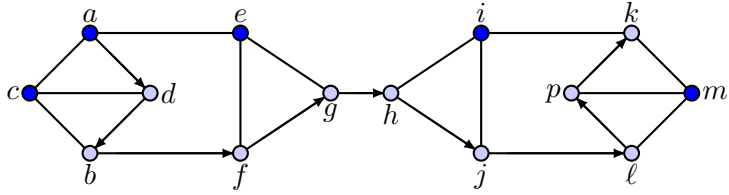

We will now assume that the diamond-unit $D^*$ contains the vertex $\ell$. Thus, $G$ is a graph with order $n = 14$ shown in Figure~\ref{fig:main-thm-4}, where $V(D^*) = \{k, \ell, m, p \}$. Then, the set $\{a, c, e, i, m \}$ is a minimum zero forcing set of $G$, and the set $\{c, e, i, m\}$ is a minimum dominating set of $G$. Thus, 
\[
Z(G) = 5 < \gamma(G) + 2 = 7.
\]
Hence, we may assume that the vertex $k$ belongs to a triangle-unit in the $\Delta$-D-partition of $V(G)$, since otherwise the desired result follows. 

Consider now the connected, cubic, and claw-free graph $G'$ with order $n'$ obtained from $G$ by removing the vertices in the set $\{e, f, g, h, i, j\}$ and any edges incident with this set, and then adding the edges $ak$ and $b\ell$. 

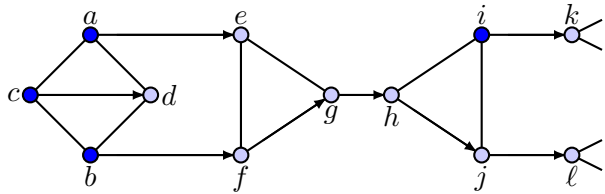
\begin{figure}[htb]
\begin{center}
\begin{tikzpicture}[scale=.8,style=thick,x=1cm,y=1cm]
\def\vr{3.5pt}

\path (1.0, 1.0) coordinate (a);
\path (0.0, 0.0) coordinate (c);
\path (1.0, -1.0) coordinate (b);
\path (2.0, 0.0) coordinate (d);

\path (3.5, 1.0) coordinate (e);
\path (3.5, -1.0) coordinate (f);
\path (5.0, 0.0) coordinate (g);

\path (7.5, 1.0) coordinate (i);
\path (7.5, -1.0) coordinate (j);
\path (6.0, 0.0) coordinate (h);

\path (9.0, 1.0) coordinate (k);
\path (9.0, -1.0) coordinate (ell);

\path (9.5, 1.25) coordinate (k1);
\path (9.5, 0.75) coordinate (k2);

\path (9.5, -0.75) coordinate (ell1);
\path (9.5, -1.25) coordinate (ell2);

\draw (a) -- (c);
\draw (a) -- (d);
\draw (c) -- (d);
\draw (c) -- (b);
\draw (d) -- (b);

\draw (a) -- (e);
\draw (b) -- (f);
\draw (e) -- (f);
\draw (e) -- (g);
\draw (f) -- (g);

\draw (g) -- (h);
\draw (h) -- (i);
\draw (h) -- (j);
\draw (i) -- (j);

\draw (i) -- (k);
\draw (j) -- (ell);

\draw (k) -- (k1);
\draw (k) -- (k2);

\draw (ell) -- (ell1);
\draw (ell) -- (ell2);

\draw[-latex, shorten >=.65mm, shorten <=3.5mm, line width=0.25mm] (c) -- (d);
\draw[-latex, shorten >=.65mm, shorten <=3.5mm, line width=0.25mm] (a) -- (e);
\draw[-latex, shorten >=.65mm, shorten <=3.5mm, line width=0.25mm] (b) -- (f);
\draw[-latex, shorten >=.65mm, shorten <=3.5mm, line width=0.25mm] (f) -- (g);
\draw[-latex, shorten >=.65mm, shorten <=3.5mm, line width=0.25mm] (g) -- (h);
\draw[-latex, shorten >=.65mm, shorten <=3.5mm, line width=0.25mm] (h) -- (j);
\draw[-latex, shorten >=.65mm, shorten <=3.5mm, line width=0.25mm] (j) -- (ell);
\draw[-latex, shorten >=.65mm, shorten <=3.5mm, line width=0.25mm] (i) -- (k);

\draw (a) [fill=blue] circle (\vr);
\draw (b) [fill=blue] circle (\vr);
\draw (c) [fill=blue] circle (\vr);
\draw (d) [fill=blue!20] circle (\vr);
\draw (e) [fill=blue!20] circle (\vr);
\draw (f) [fill=blue!20] circle (\vr);
\draw (g) [fill=blue!20] circle (\vr);
\draw (h) [fill=blue!20] circle (\vr);
\draw (i) [fill=blue] circle (\vr);
\draw (j) [fill=blue!20] circle (\vr);
\draw (k) [fill=blue!20] circle (\vr);
\draw (ell) [fill=blue!20] circle (\vr);

\node at (a) [above] {$a$};
\node at (b) [below] {$b$};
\node at (c) [left] {$c$};
\node at (d) [right] {$d$};
\node at (e) [above] {$e$};
\node at (f) [below] {$f$};
\node at (g) [below] {$g$};
\node at (h) [below] {$h$};
\node at (i) [above] {$i$};
\node at (j) [below] {$j$};
\node at (k) [above] {$k$};
\node at (ell) [below] {$\ell$};

\end{tikzpicture}
\end{center}
\caption{The subgraph of $G$ appearing towards the end of the proof given for Theorem~\ref{conj:confirmed}.}
\label{fig:hard-subcase-1}
\end{figure}

\begin{claim}\label{claim:z-claim-2}
$Z(G) \leq Z(G') + 2$
\end{claim}
\proof Let $S'$ be a minimum zero forcing set of $G'$ and let $D$ be the diamond-unit of $G'$ (and also of $G$) with vertex set $V(D) = \{a, b, c, d\}$. We will consider the following three subclaims. Suppose first that $|V(D) \cap S'| = 3$. In this case Let $S = S'\cup \{i\}$. Note that any three vertices of $D$ will zero force all vertices of $D$ to become colored blue, and so, we assume $S'\cap D = \{a, b, c\}$. Then, as shown in Figure~\ref{fig:hard-subcase-1}, all vertices of the subgraph $H$ in $G$ become colored blue. Then, after $k$ and $j$ become colored blue, all color changes in $G'$ will happen in $G$. Thus, $S$ is a zero forcing set of $G$. Hence, $Z(G) \leq |S| = |S'| + 1 = Z(G') + 1 < Z(G') + 2 $. 

Next suppose $|V(D) \cap S'| = 2$. In this case, we may assume without loss of generality that $a \notin S'$ and let $S = S'\cup \{a, i\}$. Then, as shown in Figure~\ref{fig:hard-subcase-1}, all vertices of the subgraph $H$ in $G$ become colored blue. Then, after $k$ and $j$ become colored blue, all color changes in $G'$ will happen in $G$. Thus, $S$ is a zero forcing set of $G$. Hence, $Z(G) \leq |S| = |S'| + 2 \leq Z(G') + 2 $. 

\begin{figure}[htb]
\begin{center}
\begin{tikzpicture}[scale=.8,style=thick,x=1cm,y=1cm]
\def\vr{3.5pt}

\path (1.0, 1.0) coordinate (a);
\path (0.0, 0.0) coordinate (c);
\path (1.0, -1.0) coordinate (b);
\path (2.0, 0.0) coordinate (d);

\path (3.5, 1.0) coordinate (e);
\path (3.5, -1.0) coordinate (f);
\path (5.0, 0.0) coordinate (g);

\path (7.5, 1.0) coordinate (i);
\path (7.5, -1.0) coordinate (j);
\path (6.0, 0.0) coordinate (h);

\path (9.0, 1.0) coordinate (k);
\path (9.0, -1.0) coordinate (ell);

\path (9.5, 1.25) coordinate (k1);
\path (9.5, 0.75) coordinate (k2);

\path (9.5, -0.75) coordinate (ell1);
\path (9.5, -1.25) coordinate (ell2);

\draw (a) -- (c);
\draw (a) -- (d);
\draw (c) -- (d);
\draw (c) -- (b);
\draw (d) -- (b);

\draw (a) -- (e);
\draw (b) -- (f);
\draw (e) -- (f);
\draw (e) -- (g);
\draw (f) -- (g);

\draw (g) -- (h);
\draw (h) -- (i);
\draw (h) -- (j);
\draw (i) -- (j);

\draw (i) -- (k);
\draw (j) -- (ell);

\draw (k) -- (k1);
\draw (k) -- (k2);

\draw (ell) -- (ell1);
\draw (ell) -- (ell2);

\draw[-latex, shorten >=.65mm, shorten <=3.5mm, line width=0.25mm] (ell) -- (j);
\draw[-latex, shorten >=.65mm, shorten <=3.5mm, line width=0.25mm] (j) -- (h);
\draw[-latex, shorten >=.65mm, shorten <=3.5mm, line width=0.25mm] (i) -- (k);
\draw[-latex, shorten >=.65mm, shorten <=3.5mm, line width=0.25mm] (h) -- (g);
\draw[-latex, shorten >=.65mm, shorten <=3.5mm, line width=0.25mm] (g) -- (f);
\draw[-latex, shorten >=.65mm, shorten <=3.5mm, line width=0.25mm] (f) -- (b);
\draw[-latex, shorten >=.65mm, shorten <=3.5mm, line width=0.25mm] (b) -- (d);
\draw[-latex, shorten >=.65mm, shorten <=3.5mm, line width=0.25mm] (e) -- (a);

\draw (a) [fill=blue!20] circle (\vr);
\draw (b) [fill=blue!20] circle (\vr);
\draw (c) [fill=blue] circle (\vr);
\draw (d) [fill=blue!20] circle (\vr);
\draw (e) [fill=blue] circle (\vr);
\draw (f) [fill=blue!20] circle (\vr);
\draw (g) [fill=blue!20] circle (\vr);
\draw (h) [fill=blue!20] circle (\vr);
\draw (i) [fill=blue] circle (\vr);
\draw (j) [fill=blue!20] circle (\vr);
\draw (k) [fill=blue!20] circle (\vr);
\draw (ell) [fill=blue!20] circle (\vr);

\node at (a) [above] {$a$};
\node at (b) [below] {$b$};
\node at (c) [left] {$c$};
\node at (d) [right] {$d$};
\node at (e) [above] {$e$};
\node at (f) [below] {$f$};
\node at (g) [below] {$g$};
\node at (h) [below] {$h$};
\node at (i) [above] {$i$};
\node at (j) [below] {$j$};
\node at (k) [above] {$k$};
\node at (ell) [below] {$\ell$};

\end{tikzpicture}
\end{center}
\caption{The subgraph of $G$ appearing towards the end of the proof given for Theorem~\ref{conj:confirmed}.}
\label{fig:hard-subcase-2}
\end{figure}
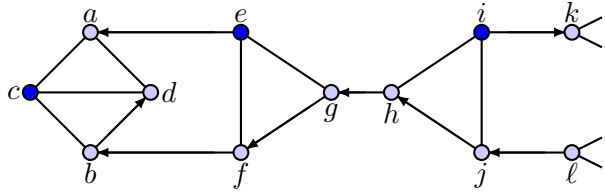

Suppose $|V(D) \cap S'| = 1$, and so, either $c$ or $d$ belong to $S'$. If Moreover, either $a$ is forced to become colored blue by $i$ in $G'$, or $b$ is forced to become colored blue by $\ell$ in $G'$. Without loss of generality suppose $\ell$ forces $b$ to become colored blue in $G'$. Under this assumption we let $S = S' \cup \{e, i\}$. Then in $G$, at some point of the zero forcing process starting with with $S$ as our initial set of blue colored vertices, $\ell$ would force $j$ to become colored blue, thereafter, all vertices of $H$ become colored blue as shown in Figure~\ref{fig:hard-subcase-2}. \smallqed 
 
\begin{claim}\label{claim:dom-claim-2}
$\gamma(G') \leq \gamma(G) - 2$
\end{claim}
\proof Observe there exists some minimum dominating set of $G$ containing the set of vertices $\{c, g, k, \ell\}$. If $X$ is this dominating set, then the let $X' = X \setminus \{i, g\}$ is a dominating set of $G'$. Thus, 
\[
\gamma(G') \leq |X'| = |X| - 2 = \gamma(G) - 2. 
\]
\smallqed 

Now that we have established our final claims, apply the inductive hypothesis to the graph $G'$, and observe
\[
Z(G') < \gamma(G') + 2.
\]
Then by Claim~\ref{claim:z-claim-2}, $Z(G) \leq Z(G') + 2$, and by Claim~\ref{claim:dom-claim-2}, $\gamma(G') = \gamma(G) -2$. Thus, 
\[
Z(G) - 2= Z(G') < \gamma(G') + 2 \leq (\gamma(G) - 2) + 2.
\]
Hence,
\[
Z(G) < \gamma(G) + 2,
\]
and the proof of the desired result is finished. \qed


\section{Conclusion}
In this paper, we have proven $Z(G) \leq \gamma(G) + 2$ for any connected, cubic, and claw-free graph $G$. This result resolves Conjecture~\ref{conj:confirmed} in the affirmative and so provides further support for the usefulness of the artificial intelligence program TxGraffiti. Moreover, since $\gamma(G) \leq \alpha(G)$ for all graphs, Theorem~\ref{thm:main} also improves on the bound $Z(G)\leq \alpha(G) + 1$ for connected, cubic, and claw-free graphs given in~\cite{Davila2, DaHe19c}. Indeed, 
\[
Z(N_k) = \gamma(N_k) + 2 < \alpha(N_k) + 1,
\]
for all diamond-necklaces $N_k$ with $k \geq 3$. Furthermore, when the additional condition that $G$ is claw-free is imposed, Theorem~\ref{thm:main} also improves on the more general bound $Z(G) \leq 2\gamma(G)$ for any connected and cubic graph $G$, other than $K_4$, presented in~\cite{Davila2, DaHe21a}.

By Theorem~\ref{thm:main}, if $G$ is connected, cubic, and claw-free, then $Z(G) = \gamma(G) + 2$, if and only if $G \in \mathcal{N}_{\text{cubic}}^*$. 
It therefore seems interesting to find all connected, cubic, and claw-free graphs for which $Z(G) = \gamma(G) + 1$. 
\begin{quest}\label{quest:1}
What connected, cubic, and claw-free graphs with $G\notin \mathcal{N}_{\text{cubic}}^*$ satisfy $Z(G) \leq \gamma(G) + 1$ with equality?
\end{quest}

Surprisingly, after completing the proof of Theorem~\ref{thm:main}, we discovered the following conjecture of TxGraffiti while prompting it to directly conjecture on the zero forcing number for all possible combinations of the hypothesis that included cubic graphs.  
\begin{conj}[TxGraffiti -- Open]\label{conj:new}
If $G$ is a connected, cubic, and diamond-free graph, then
\[
Z(G) \leq \gamma(G) + 2,
\]
and this bound is sharp. 
\end{conj}
Conjecture~\ref{conj:new} was ranked as more substantial than the more well-known $\alpha$-Z Conjecture, so this conjecture also warrants further investigation.

\end{document}